\documentclass{amsart}

\usepackage{graphicx}
\usepackage{multirow}
\usepackage{multicol}
\usepackage{graphicx}
\usepackage{subcaption}
\usepackage[bookmarks=false]{hyperref}

\usepackage{enumerate}
\usepackage{enumitem}
\usepackage{longtable}
\usepackage{geometry}
\usepackage{tabularx}
\usepackage{lipsum}
\usepackage[linesnumbered,boxed]{algorithm2e}

\usepackage{mdframed}
\usepackage{xcolor}
\usepackage{graphicx}
\usepackage{amsmath}
\usepackage{lscape} 

\usepackage{fixltx2e}
\usepackage{bm}
\usepackage{amsthm}
\usepackage{thmtools}
\usepackage{float}
\usepackage{slashbox}
\usepackage{comment}
\usepackage{tabularx} 
\usepackage{enumitem}
\usepackage{amsmath}
\usepackage{placeins}

\DeclareMathOperator{\len}{length}
\DeclareMathOperator{\isp}{isprime}
\SetKwProg{Fn}{Function}{}{end}\SetKwFunction{Prep}{Prepare}
\SetKwProg{Fn}{Function}{}{end}\SetKwFunction{Erato}{SmallPrimes}
\SetKwProg{Fn}{Function}{}{end}\SetKwFunction{Gmqr}{Generatem2qr}
\SetKwProg{Fn}{Function}{}{end}\SetKwFunction{GGC}{GGC}
\SetKwProg{Fn}{Function}{}{end}\SetKwFunction{Checkone}{Check1}
\SetKwProg{Fn}{Function}{}{end}\SetKwFunction{Checktwo}{Check2}

\SetKwProg{Fn}{Function}{}{end}\SetKwFunction{GGCone}{GGC1}
\SetKwProg{Fn}{Function}{}{end}\SetKwFunction{GGCtwo}{GGC2}
\SetKwProg{Fn}{Function}{}{end}\SetKwFunction{Genismp}{GenerateIsm1p}
\SetKwProg{Fn}{Function}{}{end}\SetKwFunction{Genresp}{GenerateResiduePattern}

\SetKwProg{Fn}{Function}{}{end}\SetKwFunction{Gmpr}{Generatem1pr}
\SetKwProg{Fn}{Function}{}{end}\SetKwFunction{Gismq}{Generateism2q}

\SetKwInOut{Input}{Input}
\SetKwInOut{Output}{Output}
\SetAlgoLongEnd

\newtheorem{theorem}{Theorem}[section]
\newtheorem{lemma}[theorem]{Lemma}
\newtheorem{conjecture}[theorem]{Conjecture}
\newtheorem{proposition}[theorem]{Proposition}
\newtheorem{corollary}[theorem]{Corollary}

\theoremstyle{definition}
\newtheorem{definition}[theorem]{Definition}
\newtheorem{example}[theorem]{Example}

\theoremstyle{remark}
\newtheorem{remark}[theorem]{Remark}

\newcommand{\cor}{\begin{corollary}}
	\newcommand{\roc}{\end{corollary}}
\newcommand{\pr}{\begin{proof}}
	\newcommand{\rp}{\end{proof}}
\newcommand{\ex}{\begin{example}}
	\newcommand{\xe}{\end{example}}
\newcommand{\tm}{\begin{theorem}}
	\newcommand{\mt}{\end{theorem}}
\newcommand{\lm}{\begin{lemma}}
	\newcommand{\ml}{\end{lemma}}
\newcommand{\df}{\begin{definition}}
	\newcommand{\fd}{\end{definition}}
\newcommand{\prop}{\begin{proposition}}
	\newcommand{\porp}{\end{proposition}}
\newcommand{\cj}{\begin{conjecture}}
	\newcommand{\jc}{\end{conjecture}}
\newcommand{\rmk}{\begin{remark}}
	\newcommand{\kmr}{\end{remark}}

\theoremstyle{plain}
\newtheorem*{GGC1}{Generalised Goldbach Conjecture (GGC)}

\theoremstyle{plain}
\newtheorem*{RGGC}{Reduced Form of the Generalised Goldbach Conjecture (RGGC)}

\theoremstyle{plain}
\newtheorem*{EGGC}{Extension of the Generalised Goldbach Conjecture (EGGC)}

\numberwithin{equation}{section}

\begin{document}

\title[Empirical verification of a new generalisation of Goldbach's conjecture]{Empirical verification of a new generalisation of Goldbach's conjecture up to $10^{12}$ (or   $10^{13}$)  for all coefficients  $\leq 40$} 
	


\author{Zsófia Juhász}
\address{Eötvös Loránd University, Dept. of Computer Algebra, Faculty of Informatics,  Budapest, Hungary}
\curraddr{}
\email{jzsofia@inf.elte.hu}
\thanks{}

\author{Máté Bartalos}
\address{}
\curraddr{}
\email{bartmate@gmail.com}
\thanks{}

\author{Péter Magyar}
\address{Eötvös Loránd University Center Savaria, Faculty of Informatics, Szombathely, Hungary}
\curraddr{}
\email{map115599@gmail.com}
\thanks{}

\author{Gábor Farkas}
\address{Eötvös Loránd University Center Savaria, Faculty of Informatics, Szombathely, Hungary}
\curraddr{}
\email{farkasg@inf.elte.hu}
\thanks{}



\date{}

\dedicatory{}

\begin{abstract}
	A new generalisation of Goldbach's conjecture ($GGC$) -- also generalising that of Lemoine -- is tested, introduced by the first author. 
	It states that for every pair of positive integers $m_1, m_2$, every sufficiently large integer $n$ satisfying certain simple criteria can be expressed as $n=m_1p+m_2q$ for some primes $p$ and $q$. $GGC$ is checked 
	up to $10^{12}d$ for all (up to $10^{13}d$ for some) pairs of coefficients $m_1, m_2$, where $d=\gcd{(m_1, m_2)}$ and $m_1/d, m_2/d \leq 40$. 
	The largest counterexamples 
	 found that cannot be obtained in this form are presented. 
	 Their relatively small sizes support the plausibility of $GGC$. Lemoine's conjecture is verified up to a new record of $10^{13}$. 
	Four naturally arising verifying algorithms are described,  and their running times compared for every 
	 $m_1\leq m_2\leq 40$ relatively prime.
	These seek to find either the $p$- 
	or the $q$-minimal $(m_1, m_2)$-partitions of all numbers tested, by either a descending or an ascending search for the prime to be maximised or minimised, respectively, in the partitions. 
	For all $m_1, m_2$  descending searches  were faster than ascending ones. A heuristic explanation is provided.
	The relative speed of ascending [descending] searches for the $p$- and for the $q$-minimal partitions, respectively,  varied by $m_1, m_2$. Using the average 
	of $p^*_{m_1, m_2}(n)$ -- the minimal $p$ in all $(m_1, m_2)$-partitions of 
	$n$ -- up to a sufficiently large threshold, two functions of $m_1, m_2$ are introduced, which may help predict these rankings and could   inform new verification efforts.  Our predictions correspond well with actual  rankings.  
	 These could potentially be further improved by developing approximations to $p^*_{m_1, m_2}(n)$. Numerical data are presented, including average and maximum values of $p^*_{m_1, m_2}(n)$ up to $10^9$.	
	 An extension of GGC is proposed, 
	  also generalising the Twin prime conjecture and the assertion that there are infinitely many Sophie Germain primes.

\end{abstract}

\maketitle


\section{Introduction}	\label{Sect_Intro}





One of the best known and longest standing 
open problems 
in number theory is 
posed by the even (or strong) Goldbach  conjecture.  
 First mentioned in 1742 by C. Goldbach in a letter to L. Euler  \cite{Fuss}, it states -- in its modern form  --  that every even number greater than $2$ can be expressed as the sum of two primes.   Search for its proof or disproof has fascinated generations of scholars and 
 curious minds since. 

Progress achieved includes W. C. Lu's showing   that the number of even integers up to  $x$ which do not have Goldbach partitions is $O(x^{0,879})$ \cite{Wen}. In \cite{Chen} J. R. Chen proved that every sufficiently large even number is the sum of a prime and a semiprime (the product of at most two primes).  In 2013  H. A.  Helfgott gave a proof for the  odd (weak or ternary) Goldbach conjecture -- a weaker statement than the even Goldbach conjecture --  claiming that every odd number greater than $5$ is the sum of three primes \cite{Helfgott2013}, \cite{Helfgott2015}.\footnote{\cite{Helfgott2013} has not been published in a peer-reviewed journal, \cite{Helfgott2015} has already been accepted for publication.}



With a general proof out of reach, several efforts have been made to verify the even Goldbach conjecture ($GC$) empirically up to increasing limits \cite{Pipping}, \cite{Richstein}, \cite{Sinisalo}, \cite{Granville1989}. The current record of $4\cdot 10^{18}$ was achieved by Olivi\'era e Silva \textit{et al.}  in a large scale computational project in 2014 \cite{Oliveira2014}. The  Goldbach partition of $n$  containing the smallest value of $p$ is called the minimal partition of $n$,  and the corresponding values of $p$ and $q$ are denoted by $p(n)$ and $q(n)$, respectively \cite{Granville1989}, \cite{Oliveira2014}. 
 In 
  \cite{Oliveira2014} verification   was carried out by segments of size $10^{12}$,  
   and in each interval the minimal Goldbach partitions of even numbers  were searched for using an efficient sieve method. Subsequently, outstanding values $n$ were handled individually   by `ascending search' for $p(n)$. For each interval to be tested primes -- potential candidates for $q$ -- in a somewhat larger interval were generated first, using a cache-efficient modified segmented sieve of Eratosthenes. 

   The rate of growth of $p(n)$ is of some theoretical interest. In \cite{Granville1989}  $p(n)=O(\log ^2 n\log \log n)$ was conjectured. A. Granville suggested two more precise, incompatible conjectures of the form $p(n)\leq (C+o(1))\log^2 n\log{\log{n}}$, where $C$ is `sharp' in the sense that $C$ is the smallest constant with this property: one with $C=C_2^{-1}\approx 1,51478$ and the other one with $C=2e^{-\gamma}C_2^{-1}\approx 1,70098$, where $C_2\approx 0,66016$ is the twin prime constant and $\gamma\approx 0,57722$ is the Euler-constant \cite{Oliveira2014}. Empirical comparison of their plausibility in \cite{Oliveira2014} was inconclusive due to the requirement of data up to even higher limits.
 

In 1894 $\mathrm{\acute{E}}$. Lemoine proposed a stronger version of the weak Goldbach conjecture   \cite{Lemoine},   stating that every odd number 
$n>5$ can be expressed as $n=p+2q$ for some primes $p$ and $q$. 
The highest treshold of  verification of Lemoine's conjecture ($LC$) the authors have found claims of 
is  $10^{10}$
\cite{MaketheBrain}.

In \cite{FarkasJuhasz} a new generalisation of the even Goldbach conjecture ($GGC$) 
was introduced, 
also generalising $LC$. 
It states that for every positive integer $m_1$ and $m_2$, every sufficiently large integer $n$ satisfying certain simple conditions can be expressed as $n=m_1p+m_2q$ for some primes $p$ and $q$. 
To the authors'  knowledge, apart from its special cases --   $GC$ when $m_1=m_2=1$ and $LC$ when $m_1=1, m_2=2$ --  $GGC$ has not been mentioned in the literature in its general form, except for \cite{FarkasJuhasz}; hence current paper is the second one investigating 
 it. 
 In \cite{FarkasJuhasz} $GGC$ was tested up to $10^{9}$ for each 
  $m_1, m_2 \leq 25$ relatively prime,  and the smallest value of $n$ satisfying the conditions of $GGC$ starting from which all integers $\leq 10^9$ also satisfying these 
  can be $(m_1, m_2)$-partitioned was provided.


We extend the scope and limit of verification of $GGC$ 
to all pairs of coefficients  
$m_1, m_2\leq 40$ up to  $ 10^{12}$ (up to $10^{13}$ for some $m_1, m_2$),   presenting  the greatest values of $n\leq 10^{12}$ satisfying the conditions of $GGC$ 
which cannot be $(m_1, m_2)$-partitioned 
\footnote{By this $GGC_{m_1, m_2}$ is also tested and the largest counterexample is determined up to $10^{12}d$ for every $m_1, m_2$ such that $m_1/d, m_2/d\leq 40$ where $d=\gcd{(m_1, m_2)}$, see Section \ref{Sect_Prelim}.}.

It is sufficient to consider 
the cases when  
$m_1$ and $m_2$ are relatively prime. The relatively small sizes of the largest counterexamples support $GGC$.
$LC$ is  confirmed up to a new record of $10^{13}$. Four different verifying algorithms with naturally arising designs were applied  to every pair $m_1< m_2$.\footnote{If $m_1=m_2$ then we have only two different approaches, hence only two different algorithms were applied when $m_1=m_2=1$.}  
We compare their speed for each $m_1, m_2$, provide heuristic explanations for their speed rankings, and seek predictions for the fastest one  when testing  up to large tresholds. In this paper we are not aiming to fully optimize our algorithms, but interested in comparing four natural approaches to testing. 
For each pair $m_1, m_2$,
 the fastest one  can be further improved  and potentially combined with other -- perhaps more efficient, e.g. sieving -- methods  
for testing up to higher limits in the future. 

 After 
 preliminaries 
 in Section \ref{Sect_Prelim},  the four  algorithms  
   are described in Section \ref{Sect_DescAlg}. 
 Searching for the minimal Goldbach partition at the verification of $GC$  \cite{Oliveira2014} has two analogues when checking  $GGC$ 
 with $m_1\neq m_2$: finding either the $p$- or the $q$-minimal $(m_1, m_2)$-partitions of numbers. In either case one can search in descending order for the prime to be maximised or in ascending order for the prime to be minimised  in the partitions. 
 These considerations yield four approaches to testing. Some findings about the functions $p^*_{m_1, m_2}(n)$ 
 and about the largest numbers $\hat{k}_{m_1, m_2}$ found satisfying the conditions of $GGC$ that cannot be  $(m_1, m_2)$-partitioned, which are relevant to the designs of the algorithms are also presented.   

Section \ref{Sect_Impl} provides information about the implementation of the algorithms and the measures taken to check the correctness of our
 computations. 

In Section \ref{Sect_Times}  the results regarding the speed ranking of the four algorithms for each pair $m_1\leq m_2\leq 40$ relatively prime 
-- presented in Section \ref{Sect_Tables} -- are discussed with some 
 heuristic explanations by the first author.
Summary data on running times is included.  
 Since  primes  among larger numbers   are scarcer on average, one may 
hypothesize that descending search for the prime to be maximised in the partition is faster than ascending search for the prime to be minimised. This 
is fully supported by our data. 
According to the results, whether descending [ascending] search for  the $p$- or for the $q$-minimal partitions is faster depends on the pair $m_1, m_2$. Two hypotheses using two functions of $m_1, m_2$ and of the average of $p^*_{m_1, m_2}(n)$ taken up to a sufficiently large treshold are proposed to predict these rankings. Predicted and actual rankings are compared, revealing reasonably good match.  Approximations for the functions $p^*_{m_1, m_2}(n)$ would be required for estimating  the time complexities of the  algorithms, and would hence 
help ascertain the plausibility of the hypotheses. 




In Section \ref{Sect_FGGC}  an extension of $GGC$ is suggested by the first author. Section \ref{Sect_Future} outlines our conclusions and some questions for future work.

Section \ref{Sect_Tables}  contains a subset of the data generated. 
The largest value $n\leq 10^{12}$ satisfying the conditions of $GGC$ 
that cannot be $(m_1, m_2)$-partitioned are presented for every $m_1, m_2\leq 40$, and the maximum and average values of $p^*_{m_1, m_2}(n)$ when $n\leq 10^9$ for every $m_1, m_2\leq 20$ relatively prime. 
Actual speed rankings  of the four algorithms and the speed rankings predicted by our hypotheses  are shown for every $m_1<m_2\leq 40$ relatively prime. 


Pseudocodes of the algorithms 
are attached 
in the Appendix.

\section{Preliminaries} \label{Sect_Prelim}




For any integers $a$ and $b$, $\gcd{(a, b)}$ shall denote the greatest common divisor of $a$ and $b$. The following conjecture was introduced in \cite{FarkasJuhasz}:

\begin{GGC1}
 \label{con_GGC}
	Let $m_1$ and $m_2$ be positive integers.  Then for every sufficiently large integer $n$ satisfying the conditions:
	
	\begin{enumerate}
		\item $\gcd(n, m_1)=\gcd(n, m_2) =\gcd(m_1, m_2)$ and   \label{con_GGC_1}
		\item $n\equiv m_1+m_2$   $\pmod{2^{s+1}}$, where $2^s$ is the largest power of $2$ that is a common divisor of $m_1$ and $m_2$,  \label{con_GGC_2}
	\end{enumerate}
	there exist primes $p$ and $q$ such that: 
	\begin{equation} \label{Eqn_GGpartition}
	n=m_1p+m_2q.
	\end{equation}
\end{GGC1}


The claim of $GGC$ for a given pair of coefficients $m_1, m_2$ shall be denoted by $GGC_{m_1, m_2}$. Note that $GGC_{1, 1}$ and $GGC_{1, 2}$ are Goldbach's and Lemoine's conjectures, respectively.





\df
An expression of the form \ref{Eqn_GGpartition} where $p$ and $q$ are primes
is called an $(m_1, m_2)$\textit{-Goldbach partition} (or $(m_1, m_2)$\textit{-partition}) \textit{of} $n$.   We say that 
  $n$ \textit{can be} $(m_1, m_2)$\textit{-partitioned} if it possesses at least one $(m_1, m_2)$-partition.
\fd

For any $m_1, m_2$, $n=m_1+m_2$  satisfies the conditions of $GGC_{m_1, m_2}$ 
and cannot be $(m_1, m_2)$-partitioned. Hence, if $GGC_{m_1, m_2}$ is true 
then there exists a largest positive integer satisfying the conditions of $GGC_{m_1, m_2}$  that  cannot be $(m_1, m_2)$-partitioned, which we denote by $k_{m_1, m_2}$. While $\hat{k}_{m_1, m_2}$ shall stand for the 
largest integer $\leq 10^{12}$ satisfying the conditions of $GGC_{m_1, m_2}$ that cannot be $(m_1, m_2)$-partitioned. 
We conjecture that $\hat{k}_{m_1, m_2}=k_{m_1, m_2}$ for every pair $m_1, m_2$ tested.



\df
If  $n$ can be $(m_1, m_2)$-partitioned then the smallest and the largest values of $p$ [$q$] in all $(m_1, m_2)$-partitions of $n$ are denoted by $p^*_{m_1, m_2}(n)$ [$q^*_{m_1, m_2}(n)$] and $p^{**}_{m_1, m_2}(n)$  [$q^{**}_{m_1, m_2}(n)$], respectively. We call $n=m_1p^*_{m_1, m_2}(n)+m_2q^{**}_{m_1, m_2}(n)$ the $p$\textit{-minimal}  (or $q$\textit{-maximal}) and  
$n=m_1p^{**}_{m_1, m_2}(n)+m_2q^{*}_{m_1, m_2}(n)$ the $p$\textit{-maximal}  (or $q$\textit{-minimal}) $(m_1, m_2)$\textit{-partition of} $n$.

\fd

Clearly, for any $m_1$, $m_2$  the conditions  of $GGC_{m_1, m_2}$ and $GGC_{m_2, m_1}$ on $n$ are equivalent, and every $(m_1, m_2)$-partition of $n$ is also  an $(m_2, m_1)$-partition if the order of terms is disregarded. Hence $n$  can be $(m_1, m_2)$-partitioned if and only if it can be $(m_2, m_1)$-partitioned, and in this case $p^*_{m_1, m_2}(n)=q^*_{m_2, m_1}(n)$ and $p^{**}_{m_1, m_2}(n)=q^{**}_{m_2, m_1}(n)$. Also, $\hat{k}_{m_1, m_2}=\hat{k}_{m_2, m_1}$, $GGC_{m_1, m_2}$ and $GGC_{m_2, m_1}$ are equivalent, and if they hold then $k_{m_1, m_2}=k_{m_2, m_1}$.




\prop \label{Prop_PropertiesGGC}
Let $n, m_1, m_2$ and $d$ be positive integers. Then:

\begin{enumerate}
	\item $n$ satisfies the conditions of $GGC_{m_1, m_2}$ if and only if $n'=dn$ satisfies the conditions of $GGC_{dm_1, dm_2}$, \label{Prop_PropertiesGGC_1}
	\item $n$ can be $(m_1, m_2)$-partitioned if and only if  $n'=dn$ can be $(dm_1, dm_2)$-partitioned, and in this case $p^*_{dm_1, dm_2}(n')=p^*_{m_1, m_2}(n)$ and $q^{**}_{dm_1, dm_2}(n')=q^{**}_{m_1, m_2}(n)$ and   \label{Prop_PropertiesGGC_2}
	\item  $GGC_{m_1, m_2}$ is true if and only if $GGC_{dm_1, dm_2}$ is, and in this case $k_{dm_1, dm_2}=dk_{m_1, m_2}$.  \label{Prop_PropertiesGGC_3}
	
\end{enumerate}



\porp

\pr \ 
\begin{enumerate}
	\item {
		Cearly, $\gcd (m_1, m_2)=\gcd (n, m_1)=\gcd (n, m_2) \Leftrightarrow \gcd (dm_1, dm_2)=\gcd (dn, dm_1)=\gcd (dn, dm_2)$. Let $2^s$ be the greatest power of $2$ which is a common divisor of $m_1$ and $m_2$, and $2^t$ be the greatest power of $2$ which is a divisor of $d$. Then the greatest power of $2$ which is a common divisor of $dm_1$ and $dm_2$ is $2^{s+t}$, and $\gcd (d, 2^{s+t+1})=2^t$, hence:
		
		$$dn \equiv  dm_1+dm_2 \pmod{2^{s+t+1}} \Longleftrightarrow n \equiv m_1+m_2 
		\pmod{2^{s+t+1}/\gcd (d, 2^{s+t+1})}		
		\Longleftrightarrow$$  $$n \equiv m_1+m_2 \pmod{2^{s+1}}. $$}
	
	\item For any primes $p$ and $q$: $n=m_1p+m_2q \Leftrightarrow dn=dm_1p+gm_2q$, hence the statement follows.
	
	\item It follows from statements \ref{Prop_PropertiesGGC_1} and \ref{Prop_PropertiesGGC_2}.
	
\end{enumerate}
\rp

For verification purposes, it is helpful to rewrite $GGC$ in the `reduced', equivalent form below:

\begin{RGGC} \label{con_RGGC}
	Let $m_1$ and $m_2$ be positive integers such that $\gcd(m_1, m_2)=1$. Then for every sufficiently large integer $n$ satisfying the conditions:
	\begin{enumerate}
		\item $\gcd(n, m_1)=\gcd(n, m_2)=1$ and   \label{con_RGGC_1}
		\item $n\equiv m_1+m_2$   ($\mathrm{mod}$   $2$),   \label{con_RGGC_2}
	\end{enumerate}
	there exist primes $p$ and $q$ such that: 
	
	$$n=m_1p+m_2q.$$
\end{RGGC}

The claim by $RGGC$ for a given pair of coefficients $m_1, m_2$ shall be denoted by $RGGC_{m_1, m_2}$. 
For  $m_1$ and $m_2$ relatively prime Conditions \ref{con_RGGC_1} and \ref{con_RGGC_2} of $RGGC_{m_1, m_2}$  are equivalent to Conditions \ref{con_GGC_1} and \ref{con_GGC_2} of $GGC_{m_1, m_2}$, respectively. By Proposition \ref{Prop_PropertiesGGC}:


\cor \label{cor_RGGCSuff}
For any positive integers $m_1, m_2$ with $\gcd{(m_1, m_2)}=d$, $GGC_{m_1, m_2}$ is true if and only if 
$RGGC_{m_1/d, m_2/d}$ is true, and in this case we have $k_{m_1, m_2}=dk_{m_1/d, m_2/d}$.  
\roc

Therefore in the study and verification of $GGC$ it is sufficient to consider the statements $RGGC_{m_1, m_2}$ where $m_1\leq m_2$  are relatively prime.


\subsection{Notations}
In the sequel  $p_i$ denotes the $i^{th}$ prime number ($i\in \mathbb{N}^+$), e.g. $p_1=2, p_2=3, $ etc.; $m_1, m_2$ and $n$ are positive integers, except for Section \ref{Sect_FGGC}, where they are not always positive. 
For any $n$, $\varphi(n)$ is the value of Euler's totient function at $n$, i.e. the number of positive integers less than or equal to $n$ that are relatively prime to $n$. For given  $m_1$ and $m_2$,  
 $lcm_{m_1, m_2}$ is the least common multiple of  $m_1, m_2$ and $2$. 
For any $L> \hat{k}_{m_1, m_2}$ such that there is at least one  $n$ satisfying the conditions of $GGC_{m_1, m_2}$ such that $\hat{k}_{m_1, m_2} <n\leq L$,  the average and the maximum values of $p^*_{m_1, m_2}(n)$ over all $\hat{k}_{m_1, m_2}< n\leq L$ satisfying the conditions of $GGC_{m_1, m_2}$ 
shall be referred to more succinctly as the \textit{average} and \textit{maximum}, respectively, \textit{of} $p^*_{m_1, m_2}$ \textit{up to} $L$. 
For any integers $a$ and $m\neq 0$, $a \mod m$ is the modulo $m$ residue of $a$.


\section{Verifying algorithms} 
\label{Sect_DescAlg}

In this section the four algorithms are described which were applied for checking $GGC_{m_1, m_2}$ up to 
$N_{m_1, m_2}\approx 10^{12}$ for every pair 
$m_1\leq m_2\leq 40$ relatively prime. (This means $490$ different pairs $m_1\leq m_2$.) Some results about the functions $p^*_{m_1, m_2}(n)$ and the values $\hat{k}_{m_1, m_2}$ are also presented.  

\subsection{Overview of the algorithms}
\subsubsection{Input and output} All algorithms verify $GGC_{m_1, m_2}$ in a segmented fashion. 
The input are  $m_1$ and $m_2$ relatively prime, the treshold of verification $N$, the length $\triangle$ of the segments to be checked at a time,  and a further, implementation dependent parameter $\alpha$. These can be set as required\footnote{Subject to the constrains  on the input provided in the outline of the algorithms.}, giving flexibility 
to our codes. In our implementation $N$   was  chosen to be the smallest multiple of $2m_1 m_2$ greater than or equal to $10^{12}$ -- denoted by $N_{m_1, m_2}$ -- and $\triangle$ 
to be the smallest multiple of $2m_1m_2$ greater than or equal to $5\cdot 10^7$.\footnote{Assuming $N$ and $\triangle$ are divisible by $2m_1 m_2$ slightly simplified our code at parts.} 
For every $n$ satisfying the conditions of $GGC_{m_1, m_2}$ the algorithms only check  if $n$ has an $(m_1, m_2)$-partition $n=m_1p+m_2q$ such that $m_1p\leq \alpha$ (or $m_2q\leq \alpha$). The output is the array \textit{residual} containing those $n\leq N_{m_1, m_2}$ satisfying the conditions of $GGC_{m_1, m_2}$ which do not possess such a partition. After an algorithm has finished, it remains to check by another method if numbers in \textit{residual} can be $(m_1, m_2)$-partitioned.

 \subsubsection{Functions $p^*_{m_1, m_2}(n)$ and the choice of $\alpha$} 
\begin{figure}[ht]
	\centering
	\begin{subfigure}{0.3\textwidth}
		\includegraphics[width=\textwidth]{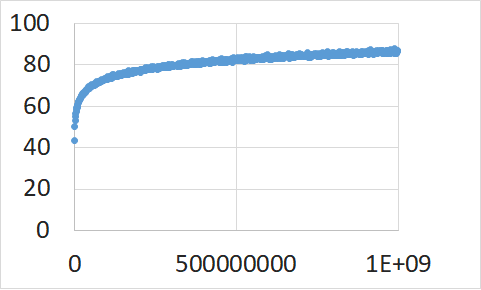}
		\caption{$m_1=1, m_2=2$}
		\label{fig:Avgp1}
	\end{subfigure}
	\hfill
	\begin{subfigure}{0.3\textwidth}
		\includegraphics[width=\textwidth]{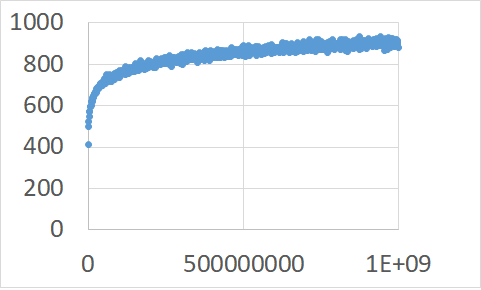}
		\caption{$m_1=4, m_2=17$}
		\label{fig:Avgp2}
	\end{subfigure}
	\hfill
	\begin{subfigure}{0.3\textwidth}
		\includegraphics[width=\textwidth]{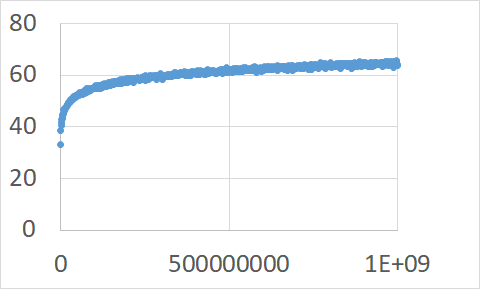}
		\caption{$m_1=7, m_2=3$}
		\label{fig:Avgp3}
	\end{subfigure}	
	\caption{The average value of $p^*_{{m_1, m_2}}(n)$ in the interval of length $10^6$ centered at $x=10^{6} k+5\cdot 10^5$ for $0\leq k\leq 10^3-1$, in cases of  $m_1, m_2$ indicated under each subfigure.}
	\label{fig:Avgp}
\end{figure}

We aimed to set the value of  $\alpha$ so that \textit{residual}  only contains numbers that cannot be $(m_1, m_2)$-partitioned at all, 
by ensuring that $m_1p^*_{m_1, m_2}(n)\leq \alpha$  for every $m_1, m_2\leq 40$ relatively prime  and $n\leq N_{m_1, m_2}$ satisfying the conditions of $GGC_{m_1, m_2}$ that can be $(m_1, m_2)$-partitioned. It was observed that  
$p^*_{m_1, m_2}(n)$ 
remains relatively small even for  large values of $n$. 
 For example, Figure \ref{fig:Avgp} demonstrates the slow growth of $p^*_{m_1, m_2}(n)$ by showing the  average of $p^*_{m_1, m_2}(n)$ in each interval of length $10^6$ centered at $x=10^6k+5\cdot 10^5$ ($0\leq k\leq 10^3-1$) in the cases $m_1=1, m_2=2$ (Subfigure \ref{fig:Avgp1}), $m_1=4 , m_2=17$ (Subfigure  \ref{fig:Avgp2}) and $m_1=7, m_2=3$ (Subfigure  \ref{fig:Avgp3}). Table \ref{Table_AverMaxp} contains the maximum and average values of $p^*_{m_1, m_2}(n)$  up to $n\leq 10^9$ for each 
$m_1, m_2\leq 20$ relatively prime. For $n\leq 10^9$, 
over all $m_1, m_2\leq 40$ relatively prime the maximum of $p^*_{m_1, m_2}(n)$ is $78697$ achieved when $m_1=32$, $m_2=37$, and the maximum of $m_1p^*_{m_1, m_2}(n)$ is $2858879$ occuring when $m_1=37$ and $m_2=38$.
Experimentally it was also found that $m_1p^*_{m_1, m_2}(n)\leq 5\cdot 10^7$ for all $n\leq N_{m_1, m_2}$ satisfying the conditions of $GGC_{m_1, m_2}$ that can be $(m_1, m_2)$-partitioned, for all $m_1, m_2\leq 40$ relatively prime. Hence in our implementation $\alpha =5\cdot 10^7$, and so   for every $m_1, m_2$, $\hat{k}_{m_1, m_2}$ is the largest number in \textit{residual}. Choosing smaller suitable $\alpha$ could have been possible, but the resulting improvements  in running times  would have been insignificant. 

\subsubsection{Values of $\hat{k}_{m_1, m_2}$}
 The value $\hat{k}_{m_1, m_2}$ for every 
$m_1, m_2\leq 40$ relatively prime  is shown in Table \ref{Table_CountEx}. 
The maximum and  average of $\hat{k}_{m_1, m_2}$ are $412987$ (reached when $m_1=32$, $m_2=37$) and  $52004,84$, respectively. 
The relatively small sizes of $\hat{k}_{m_1, m_2}$ support $GGC$. It also meant that the 
extra time required by checking numbers in \textit{residual} 
was negligible.


\subsubsection{Different approaches of the four algorithms to testing} The main difference between Algorithms 1/a, 1/b, 2/a and 2/b 
lies in their methods for checking if a number can be $(m_1, m_2)$-partitioned.  
These 
-- for given ordered pair $(m_1, m_2)$ -- 
are summarised below: 
 

\medskip








\noindent\textit{Algorithm 1/a [1/b]:} `Descending search for the prime to be maximised' in the 
partitions. 
Algorithm 1/a  [1/b] searches for the $p$-minimal [$q$-minimal] ($m_1, m_2$)-partition  $n=m_1p^*_{m_1, m_2}(n)+m_2q^{**}_{m_1, m_2}(n)$  [$n=m_1p^{**}_{m_1, m_2}(n)+m_2q^{*}_{m_1, m_2}(n)$] by trying all possible candidates $q$ [$p$] for $q^{**}_{m_1, m_2}(n)$   [for $p^{**}_{m_1, m_2}(n)$] in decreasing  order until it finds that $n-m_2q=m_1p$ [$n-m_1p=m_2q$] 
for some prime $p$ [$q$].  

\bigskip

\noindent\textit{Algorithm 2/a [2/b]:} `Ascending search for the prime to be minimised' in the 
partitions. 
Algorithm 2/a [2/b] searches for the $p$-minimal [$q$-minimal] ($m_1, m_2$)-partition  $n=m_1p^*_{m_1, m_2}(n)+m_2q^{**}_{m_1, m_2}(n)$  [$n=m_1p^{**}_{m_1, m_2}(n)+m_2q^{*}_{m_1, m_2}(n)$] by trying all possible candidates $p$ [$q$] for $p^*_{m_1, m_2}(n)$   [for $q^{*}_{m_1, m_2}(n)$] in increasing order until it finds that $n-m_1p=m_2q$ [$n-m_2q=m_1p$] 
for some prime $q$ [$p$].

\medskip


  
 



Algorithms 1/a and 1/b  [2/a and 2/b] can be implemented by the same program by interchanging the values of $m_1$ and $m_2$. Hence only Algorithms 1/a and 2/a are described in this section, referred to as Algorithms 1 and 2, respectively. 


\subsubsection{Simplified outlines of Algorithms 1 and 2} \label{Subsect_Outline}\  

\ 

\noindent \textit{Input:} 
$m_1, m_2, N,$ 
$\triangle, \alpha\in \mathbb{N}^+$ such that $\gcd (m_1, m_2)=1$, $N >9$, 
$2m_1m_2|N$, 
$2m_1m_2|\triangle$ and $\alpha\leq \triangle$. 


\medskip

\noindent \textit{Output:} array \textit{residual} containing all numbers $n \leq N$ 
satisfying the conditions of $GGC_{m_1, m_2}$ for which there are  no primes $p$ and $q$ such that $n=m_1p+m_2q$ and $m_1p\leq \alpha$.

\medskip
\begin{enumerate}
	
	\item \label{AlgI} {Phase I: Unsegmented phase
		
		\begin{enumerate}
			
			\item \label{AlgI/1/1} Generating `small' primes 
			up to  
			$K= \max{ \{\lfloor\sqrt{\smash[b]{N/m_2}}\rfloor, \lfloor\alpha/m_1 \rfloor \}}$. 
			
			\item \label{AlgI/1/2} Generating all numbers 
			$m_1p\leq \alpha$ where $p$ is prime. 
			In  Algorithm 2 these are sorted and stored separately according to their modulo $m_2$ residues. 
			

			\item \label{AlgI/2} Generating the modulo $lcm_{m_1, m_2}$ `residue wheel', i.e. the array of all 
			modulo $lcm_{m_1, m_2}$ residues relatively prime to $m_1m_2$ and congruent to   $m_1+m_2$ modulo $2$. 
			
	\end{enumerate} }

	\item \label{AlgII} {Phase II: Checking $GGC_{m_1, m_2}$ segment by segment

		For each interval $[A,B)$:

		\begin{enumerate}
			
			\item \label{AlgII/1} {Generating `large' primes 
				and their $m_2$-times multiples in an interval.  

				\begin{enumerate}
					\item \label{AlgII/1/1} Generating all primes in interval $[C/m_2, D/m_2)$. (The values $C$ and $D$ depend on $A$ and $B$.) 
					\item \label{AlgII/1/2} Generating all numbers of the form $m_2q$ in interval $[C, D)$, where $q$ is prime. In  Algorithm 1 these are sorted and stored separately according to their modulo $m_1$ residues.
					
				\end{enumerate}

				
				
			}
			\item \label{AlgII/2} Checking $GGC_{m_1, m_2}$ in interval $[A,B)$.
		\end{enumerate}
	}
	
\end{enumerate}

\bigskip

\subsubsection{Some ideas applied in both algorithms} 
For checking if every number in an interval $[A, B)$ satisfying the conditions of $GGC_{m_1, m_2}$ has a partition $m_1p+m_2q$ such that $m_1p\leq \alpha$, it is sufficient to possess the lists of all numbers $m_1p\leq \alpha$ where $p$ is prime,  and of all numbers $m_2q$ in interval  $[\max{\{0, A-\alpha\}},B)$ where $q$ is prime. These lists are generated in Phases I and II, respectively. Although 
methods with lower asymptotic time complexities  exist
\cite{GriesMisra}, \cite{Barstow}, 
\cite{Misra}, \cite{Bengelloun}, \cite{Prichard}, in Phases I and II, the sieve of Eratosthenes and a segmented version of this, respectively, is used  to generate primes.

The following observation helped speed up testing: If $n=m_1p+m_2q$ is an $(m_1, m_2)$-partition 
then

                  


\begin{minipage}{.5\linewidth}
	\begin{equation}
	m_2q\equiv  n\pmod{m_1} \textrm{\;\;\;\;\;\;\;\;\; and } \label{Eqn_1}
	\end{equation}
\end{minipage}%
\begin{minipage}{.5\linewidth}
	\begin{equation}
		m_1p\equiv n \pmod{m_2}.     \label{Eqn_2}
	\end{equation}
\end{minipage}

\bigskip

 
 

 
Therefore, for  any $n$, Algorithm 1 [2] in Phase II tries only those primes $q$ [$p$] as 
candidates for  $q^{**}_{m_1, m_2}(n)$ [$p^*_{m_1, m_2}(n)$] which satisfy congruence \ref{Eqn_1} [\ref{Eqn_2}], hence reducing
the number of candidates 
tested 
by  approximately  a factor of $1/\varphi(m_1)$ 
 [$1/\varphi(m_2)$]. 
In order to facilitate this, when generating numbers of the form $m_2q$ [$m_1p$]  in an interval  [up to $\alpha$]  Algorithm 1 [2] also sorts them by their modulo $m_1$ [$m_2$] residues.

\subsection{Detailed description of the steps} \label{Subsect_DescrSteps} 

The pseudocode of the main program 
{\tt{GGC1}} [{\tt{GGC2}}] implementing Algorithm 1 [2] and those of  procedures {\tt{GenerateIsm1p}},  {\tt{Generatem1pr}},  {\tt{Generatem2qr}}, {\tt{Generateism2q}}, {\tt{Check1}} and   {\tt{Check2}} described below 
can be found in the Appendix.


\medskip


\subsubsection{Phase I: Unsegmented phase}
\ 



\noindent (a) {\textit{Generating `small' primes:}
In both algorithms  a list of all `small' primes $\leq K$ is generated first   by procedure {\tt{SmallPrimes($K$)}} using the sieve of Eratosthenes,  where $3\leq K\in \mathbb{N}$ is an  implementation dependent treshold. Small primes are used for two purposes later: for the generation of all numbers $m_1p \leq\alpha$ where $p$ is prime in Phase I, and at the sieving for large primes by segments in Phase II, up to $N_{m_1, m_2}/m_2$.  
Therefore  $K\geq \max{\{\lfloor \sqrt{\smash[b]{N_{m_1, m_2}/m_2}} \rfloor, \lfloor\alpha/m_1\rfloor \}}$
must hold. 
We set $K=\max{\{\lfloor \sqrt{\smash[b]{N_{m_1, m_2}/m_2}} \rfloor, \lfloor\alpha/m_1 \rfloor \}}$
for every pair $m_1, m_2$.
In both algorithms the output are global arrays $\mathit{isprime}$ and $\mathit{primes}$, where $\mathit{isprime}$ is a boolean array of length $\lfloor (K-1)/2\rfloor$ 
such that for every $0\leq i\leq \lfloor(K-3)/2\rfloor$: 
$\mathit{isprime}[i]=1$ if and only of  $2i+3$ is a prime and $\mathit{isprime}[i]=0$ otherwise;  $\mathit{primes}$ contains the list of all primes less than or equal to $K$ in increasing order, i.e. for every $0\leq i\leq K-1$: $\mathit{primes}[i]=p_{i+1}$.
}

\bigskip



			
		
				
					
		

		
			
		

\noindent  (b) \textit{Generating the $m_1$-times multiples of `small' primes:}
In both algorithms all numbers $m_1p\leq \alpha$ are generated 
where $p$ is prime. 
In Algorithm 2 these are sorted according to their  modulo $m_2$ residues.
In Algorithm 1 the  boolean array $ism_1p$ of length $\alpha +1$ is generated by procedure {\tt{GenerateIsm1p}}($\alpha$) where for every $0\leq i\leq \alpha$: $ism_1p[i]=1$ if and only if $i=m_1p$ for some prime $p$. In Algorithm 2 for every $0\leq r<m_2$ 
the array $\mathit{m_1p[r]}$ is generated by procedure {\tt{Generatem1pr}}($\alpha$) containing all numbers $m_1p\leq \alpha$ (in increasing order)  
where $p$ is prime  and $r= m_1p  \mod{m_2}$. 

\bigskip


\noindent (c) \textit{Generating the modulo $lcm_{m_1, m_2}$ `residue wheel':} When checking $GGC_{m_1, m_2}$ only those numbers $n$ need to be tested 
which satisfy the conditions of $GGC_{m_1, m_2}$, which holds if and only if the residue $n \mod lcm_{m_1, m_2}$ satisfies these. By procedure {\tt{GenerateResiduePattern}}($m_1, m_2$) a boolean array $res$ of length $lcm_{m_1, m_2}$ is generated such that for every $0\leq i\leq lcm_{m_1, m_2}-1$, $res[i]=1$ if and only if $\gcd(i, m_1)=\gcd(i, m_2)=1$ and $i\equiv m_1+m_2 \pmod{2}$. This 
is used later 
for deciding if a certain $n$ needs to be tested.




		

\subsubsection{Phase II: Segmented phase:}

\ 






\noindent (a) \textit{Generating all numbers 
	$m_2q$ in an interval where $q$ is prime:}
%
%
\noindent
For given integers $0\leq C<D$ such that $2m_1m_2|C$ and $2m_1m_2|D$, procedure {\tt{Generatem2qr}}$(C, D)$ in Algorithm 1  generates all numbers of the form $m_2q$ where $q$ is prime, in interval $[C, D)$, and  stores each 
$m_2q$  
in array $\mathit{m_2q[r]}$ where $r= m_2q \mod m_1$  ($0\leq r < m_1$). 
For given integers $0\leq C<D$ such that $2m_2|C$ and $2m_2|D$  procedure {\tt{Generateism2q}}$(C, D)$ in Algorithm 2 outputs boolean array $\mathit{ism_2q}$ of length $D-C$ such that for every $0\leq i< D-C-1:$ $\mathit{ism_2q}[i]=1$  if and only if $C+i=m_2q$ for some prime $q$.

\bigskip




\noindent (b) \textit{Checking $GGC_{m_1, m_2}$ in an interval:} For given integers  $0\leq A<B$,  where $2m_1m_2| A$  and  $2m_1m_2  | B$, 
procedure {\tt{Check1}}$(A, B)$ in Algorithm 1 [{\tt{Check2}}$(A, B)$  in Algorithm 2]  checks for every $n$ in $[A,B)$ satisfying the conditions  of $GGC_{m_1, m_2}$  if there exist primes $p$ and $q$ such that $n=m_1p+m_2q$ and $m_1p\leq \alpha$. {\tt{Check1}}$(A, B)$  [{\tt{Check2}}$(A, B)$] looks for the $p$-minimal $(m_1, m_2)$-partition of $n$,  applying  a `descending' search for $q^{**}_{m_1, m_2}(n)$ [an `ascending' search for $p^{*}_{m_1, m_2}(n)$]. It looks for $m_2q^{**}(n)$  [$m_1p^{*}(n)$] by trying in decreasing [increasing] order the values $m_2q$ [$m_1p$] where $q$ [$p$] is prime  such that $m_2q\equiv n \pmod{m_1}$ [$m_1p\equiv n \pmod{m_2}$] -- taking these from array $m_2q\mathrm{[r]}$ [$m_1p$] where $r= n\mod{m_1}$  -- and checking if $n-m_2q$ [$n-m_1p$] is of the form $m_1p$ [$m_2q$] for some prime $p$ [$q$]. If such value $m_2q$ [$m_1p$] is found then $q^{**}(n)=q$  [$p^{*}(n)=p$] and $p^*(n)=(n-m_2q)/m_1$ [$q^{**}(n)=(n-m_1p)/m_2$]. If no such value $m_2q$ [$m_1p$] is found then $n$ is added to array \textit{residual}.  The output of both procedures is array \textit{residual} of those numbers $n$ in $[A,B)$ satisfying the conditions  of $GGC_{m_1, m_2}$  for which there exist no primes $p$ and $q$ such that $n=m_1p+m_2q$ and $m_1p\leq \alpha$. 

\ 
\subsubsection{The main programs}
\nopagebreak


Algorithm 1 [2] is implemented by the main program\\ {\tt{GGC1}}($N,m_1,m_2, \triangle$, $\alpha$)  
 [{\tt{GGC2}}($N,m_1,m_2, \triangle$, $\alpha$)]. 
Before performing {\tt{Check1}}($A, B$) [{\tt{Check2}}($A, B$)]
all numbers of the form $m_2q$ where $q$ is prime need to be obtained in interval $[\max{(0, A-\alpha)}, B)$. In order for this, in each iteration of loop 7-23 in Algorithm 1 [loop 7-18 in Algorithm 2], the numbers $m_2q$ are generated by \texttt{Generatem2qr} 
in step 20   [\texttt{Generateism2q} 
in step 12] only in interval $[A, B)$, and starting from the second iteration those in $[\max{\{0, A-\alpha\}}, A)$ are kept from the previous iteration in steps 9-19 [in step 10]  in arrays $m_2q[r]$ [in array $ism_2q\_old$]  and added. Therefore during both algorithms every number $m_2q\leq N$ where $q$ is prime is generated exactly once.

\smallskip

\smallskip



\section{Implementation and checking for correctness} \label{Sect_Impl}

Algorithms 1 and 2 were implemented in C++. For each 
$m_1\leq m_2\leq 40$ relatively prime, 
$GGC_{m_1, m_2}$ was checked up to $N_{m_1, m_2}$ by  Algorithms 1/a, 1/b, 2/a and 2/b.\footnote{For $m_1=m_2=1$, Algorithms 1/a and 1/b [2/a and 2/b] are identical, hence in this case only two different algorithms  were performed.}  
 Algorithms 1/a and 1/b [2/a and 2/b] were both  carried out by the 
program for Algorithm 1 [2], by interchanging the values of  $m_1$ and $m_2$ (with $m_1< m_2$ in Algorithms 1/a and 2/a). Each algorithm was run on one core of a 32-core 64-bit Intel Xeon Scalable processor.


For each pair  $m_1\leq m_2$ tested the arrays \textit{residual} produced by the four\footnote{Only  two different algorithms in case $m_1=m_2=1$.} algorithms were identical; 
the values $p^*_{m_1, m_2}(n)$ [$q^{*}_{m_1, m_2}(n)$]  and $q^{**}_{m_1, m_2}(n)$ [$p^{**}_{m_1, m_2}(n)$] for every $\hat{k}_{m_1, m_2}<n\leq  10^6$ satisfying the conditions of $GGC_{m_1, m_2}$ were also generated 
and found identical. 


\section{Comparing the running times of the algorithms} \label{Sect_Times}

\subsection{Experimental data on running times}

For each pair $m_1\leq m_2\leq 40$ relatively prime, Algorithms 1/a and 1/b were both   faster than Algorithms 2/a and 2/b,  the former two  significantly outperforming on average the latter. The speed rankings of Algorithms 1/a and 1/b  [2/a and 2/b] varied depending on the pair $m_1< m_2$. On average over all pairs $m_1\leq m_2$, Algorithms 1/a and 1/b [2/a and 2/b] showed very similar speed performances.  Table \ref{Table_FastestSlowest} presents the average, lowest and highest running times of each  algorithm and the pair $m_1, m_2$ where the latter occurred:

\begin{table}[htbp]
	\centering
	\caption{Lowest, highest and average running times (sec) of the algorithms up to $N_{m_1, m_2}\approx 10^{12}$ over all pairs $m_1\leq m_2\leq 40$ relatively prime.} 
	\label{Table_FastestSlowest}
	\begin{tabular}{|l|l|l|l|l|l|l|l|}
		\hline
		\multicolumn{1}{|c|}{\textbf{Algorithm}}&  \multicolumn{3}{c|}{\textbf{Lowest}} 
		&     \multicolumn{3}{c|}{\textbf{Highest}}  
		& \multicolumn{1}{c|}{\textbf{Average time (sec)}}
		
		\\
		& \multicolumn{1}{c|}{$\mathbf{m_1}$} & \multicolumn{1}{c|}{$\mathbf{m_2}$} & \multicolumn{1}{c|}{\textbf{time (sec)}} &         \multicolumn{1}{c|}{$\mathbf{m_1}$} & \multicolumn{1}{c|}{$\mathbf{m_2}$} & \multicolumn{1}{c|}{\textbf{time (sec)}} &   \\
		\hline
		\textbf{Alg. 1/a} & 7     & 30    & 22473 &        16    & 29    & 114177 &  56045
		\\
		\textbf{Alg. 1/b} & 6     & 35    & 23345 &       1     & 16    & 108614 &  54461 \\
		\textbf{Alg. 2/a} & 33    & 35    & 55742 &       31    & 32    & 293279 &  132154 \\
		\textbf{Alg. 2/b} & 35    & 39    & 57391 &        32    & 37    & 293734 &   134559\\
		\hline
	\end{tabular}%
	\label{tab:addlabel}%
\end{table}


For each pair $m_1< m_2$ tested the running times of the four algorithms ranked in one of the following four  orders from fastest to slowest: 

\medskip


	
	\textit{Group A:}  Algorithms 1/a, 1/b, 2/a, 2/b
	
	\textit{Group B:} Algorithms 1/a, 1/b, 2/b, 2/a

	\textit{Group C:} Algorithms 1/b, 1/a, 2/a, 2/b

	\textit{Group D:} Algorithms 1/b,  1/a, 2/b,  2/a
	

\medskip

Groups $A, B, C$ and $D$ contain $21$, $218$, $242$ and $8$ pairs, respectively, see Table \ref{Table_Class} in Section \ref{Sect_Tables}. 
The dominance of groups $B$ and $C$ raises the question whether the pairs in groups $A$ and $D$ would also move to one of the former when testing up to sufficiently large tresholds. In all $8$ pairs in group $D$ the running times of Algorithms 1/a and 1/b or those of 2/a and 2/b were `very close'. We ran all four algorithms for the pairs $(9, 32), (11, 29), (17, 19)$ and $(23, 29)$ in group $D$ -- and for six other pairs including $(1, 2)$ -- up to $\approx 10^{13}$ and the running times are shown in Table \ref{Table_Upto1013}. The speed rankings changed for all four pairs in group $D$.  According to this $(9, 32)$ moved to group $B$, $(11, 29)$ to group $C$ and $(17, 19)$ and $(23, 29)$ to group $A$. In the latter two cases the running times of Algorithms 1/a and 1/b were `very close' to each other, which makes it plausible that the pairs might move again to another group if testing until even higher tresholds. These results suggest that the remaining other four pairs in group $D$ may also leave this group in case of larger tresholds.



{\footnotesize{
\begin{table}[H]
	\centering
	\caption{Running times (sec) of the algorithms up to  $\approx 10^{12}$ and $\approx10^{13}$ for some $m_1, m_2$.} \label{Table_Upto1013}
	\begin{tabular}{|c|c|c|c|c|c||c|c|c|c|}
		\hline
		$\mathbf{m_1}$ & $\mathbf{m_2}$ & \multicolumn{4}{c||}{\textbf{Running times (sec) up to} 
			$\mathbf{\approx 10^{12}}$} & \multicolumn{4}{c|}{\textbf{Running times (sec) up to} 
			 $\mathbf{\approx 10^{13}}$} \\
		{\multirow{2}{*}{}} & {\multirow{2}{*}{}} & {\textbf{Alg. 1/a}} & {\textbf{Alg. 1/b}} & {\textbf{Alg. 2/a}} & {\textbf{Alg. 2/b }} & {\textbf{Alg. 1/a}} & {\textbf{Alg. 1/b}} & {\textbf{Alg. 2/a }} & {\textbf{Alg. 2/b}}\\
		
		\hline
		
		1     & 2     & 77192 & 104914 & 290478 & 182075 & 754817 & 1052855 & 2290673 & 1873002 \\
		
		1     & 3     & 42991 & 67452 & 106733 & 110383 & 423373 & 671771 & 1125928 & 1154959 \\
		
		1     & 5     & 56912 & 77743 & 129483 & 142355 & 555759 & 786325 & 1350077 & 1515172 \\
		
		1     & 7     & 63372 & 82353 & 137414 & 158389 & 627182 & 823687 & 1479262 & 1704191 \\
		
		1     & 9     & 44371 & 69002 & 101032 & 111152 & 431663 & 687392 & 1065525 & 1172360 \\
		
		1     & 11    & 69622 & 86173 & 143534 & 169134 & 745476 & 860608 & 1470691 & 1753254 \\
		
		9     & 32    & 43546 & 43514 & 132022 & 111549 & 549958 & 582308 & 1331716 & 1146279 \\
		
		11    & 29    & 74485 & 52092 & 148247 & 145263 & 744006 & 706943 & 1404459 & 1563704 \\
		
		17    & 19    & 55562 & 55052 & 143130 & 143104 & 736166 & 738674 & 1458811 & 1545744 \\
		
		23    & 29    & 80070 & 59988 & 155817 & 155277 & 800289 & 807146 & 1567925 & 1676299 \\
		\hline
	\end{tabular}%
	\label{tab:addlabel}%
\end{table}%
}}

\subsection{Estimations for the running times}

The significant parts of the computation in Algorithm 1 [2] are {\tt{Generatem2qr}} and {\tt{Check1}} [{\tt{Generateism2q}} and {\tt{Check2}}].
During 
all iterations of 
{\tt{Generatem2qr}} [{\tt{Generateism2q}}] all primes up to $N/m_2$ 
and their $m_2$-times multiples are generated, taking $O(N\log{\log{N}})$ 
\cite{Sorenson} and $\pi(N/m_2) \sim N/(m_2(\ln{N}-\ln{m_2}))=o(N\log{\log{N}})$
operations, respectively. 
Hence {\tt{Generatem2qr}} [{\tt{Generateism2q}}]  takes $O(N\log{\log{N}})$ 
time. 


In absence of approximations for the functions $p^*_{m_1, m_2}(n)$  it is difficult to estimate the number of operations performed by {\tt{Check1}}  [{\tt{Check2}}]. However, the following can be established. For given $m_1, m_2$ relatively prime 
the number of values $n\leq N_{m_1, m_2}$ tested -- i.e. of those satisfying the conditions of $GGC_{m_1, m_2}$ -- is approximately $\varphi(m_1m_2)N_{m_1, m_2}/lcm_{m_1, m_2}\approx 10^{12}\varphi(m_1m_2)/lcm_{m_1, m_2}$.










In Algorithm 1, for each $n$ tested, loop 9-30 in {\tt{Check1}} is iterated  by the number of candidates for $q^*_{m_1, m_2}(n)$ tried, which is approximately 
the number of primes $q$ in interval between $n/m_2 - m_1p^*_{m_1, m_2}(n)/m_2$ and $n/m_2$  
of length $m_1p^*_{m_1, m_2}(n)/m_2$
satisfying 
$m_2q\equiv n \pmod{m_1}$. 
This -- using  $\pi(x)-\pi(x-y)\approx y/\ln(x)$ \cite{Heath1988}  -- can be estimated as

\begin{equation} \label{Eqn_Est1}
	\begin{aligned}
		& \frac{m_1p^*_{m_1, m_2}(n)}{\varphi(m_1)m_2\ln (n/m_2)} \approx  
		& \frac{m_1p^*_{m_1, m_2}(n)}{\varphi(m_1)m_2\ln (n)}. 
	\end{aligned}
\end{equation}

In Algorithm 2 for each value $n$ tested the number of iterations of loop 9-19 in {\tt{Check2}} is equal to the number of candidates for $p^*_{m_1, m_2}(n)$ checked, which is  the number of primes $p$ up to $p^*_{m_1, m_2}(n)$ 
satisfying 
$m_1p \equiv n \pmod{m_2}$.
This can be estimated by 

\begin{equation} \label{Eqn_Est2}
\frac{\pi(p^*_{m_1, m_2}(n))}{\varphi(m_2)}\sim \frac{p^*_{m_1, m_2}(n)}{\varphi(m_2)\ln{p^*_{m_1, m_2}(n)}}.
\end{equation}

\subsection{Some heuristics}
Currently possessing no approximations for $p^*_{m_1, m_2}(n)$ and thus for the number of operations performed by {\tt{Check1}} and {\tt{Check2}}, it is unclear how the time complexities of {\tt{Generatem2qr}} and {\tt{Check1}} [{\tt{Generateism2q}} and {\tt{Check2}}] compare. In order to obtain empirical data, we ran Algorithm 1/a 
for a few (four) pairs $m_1\leq m_2$ up to the tresholds of approximately $10^6, 10^7, 10^8$ and $10^9$, and measured the times taken by {\tt{Check1}}  
and  {\tt{Generatem2qr}}. 
In one case {\tt{Check1}}  
took  around $66\%$ and in all other cases above $80\%$ (usually above $90\%$), whereas {\tt{Generatem2qr}} 
took in one case $16\%$, but in all others below $10\%$ and usually below $5\%$ of the total time. As the treshold increased,  Algorithm 1/a spent an increasing and a decreasing fraction of the total time on {\tt{Check1}}  
and on {\tt{Generatem2qr}}
, respectively. 

In the arguments below we shall assume that in Algorithm 1 [2] {\tt{Check1}} [{\tt{Check2}}] is the most time consuming part of the computation with higher time complexity than {\tt{Generatem2qr}}  [{\tt{Generateism2q}}], and hence the relative speed performances of Algorithms 1/a, 1/b, 2/a and 2/b
are determined by 
{\tt{Check1}} and {\tt{Check2}}.

\subsubsection{Comparing the running times of Algorithms 1/a and 2/a [1/b and 2/b]}
In \cite{Granville1989} it was conjectured  that $p(n)=p_{1, 1}^*(n)=O(\log ^2 n\log \log n)$, implying $p_{1, 1}^*(n)=o(n^{\varepsilon})$ for every $\varepsilon\in \mathbb{R}^+$. Based on our data we also conjecture that for every $m_1$ and $m_2$: $p^*_{m_1, m_2}(n)=o(n^{\epsilon})$ for every $\varepsilon\in \mathbb{R}^+$. Using this assumption $\ln{p^*_{m_1, m_2}(n)}=o(\ln(n))$, hence

$$ \frac{m_1p^*_{m_1, m_2}(n)}{\varphi(m_1)\ln (n)} = o\left(\frac{p^*_{m_1, m_2}(n)}{\varphi(m_2)\ln{p^*_{m_1, m_2}(n)}}\right), $$


which heuristically suggests that Algorithm 1/a [1/b] is faster than Algorithm 2/a [2/b] for every $m_1, m_2$, when run until sufficiently large treshold. This prediction is in complete accordance with our results: for each pair $m_1, m_2$ tested  Algorithms 1/a and 1/b were both faster than Algorithms 2/a and 2/b. 

\subsubsection{Comparing the running times of Algorithms 1/a and 1/b  [2/a and 2/b]}

Since for given $m_1, m_2$, in Algorithms 1/a and 1/b [2/a and 2/b]  {\tt{Check1}} [{\tt{Check2}}] checks the same number of values $n$,  one may attempt to explain their relative speed performances using some  estimate of the `average' time spent by {\tt{Check1}} [{\tt{Check2}}] on processing each value $n$. Based on estimates \ref{Eqn_Est1} and \ref{Eqn_Est2},    we introduce the following functions for any sufficienty large number $L$:

$$ f_{{}_L}(m_1, m_2):=\frac{m_1\overline{p^*}_{{}_{L}}(m_1, m_2)}{\varphi(m_1)m_2} \textrm{\hspace{2pt} and \hspace{2pt}} g_{{}_L}(m_1, m_2):= \frac{\overline{p^*}_{{}_{L}}(m_1, m_2)}{\varphi(m_2)\ln{\overline{p^*}_{{}_{L}}(m_1, m_2)}},$$

where $\overline{p^*}_{{}_{L}}(m_1, m_2)$ is the average of $p^*_{m_1, m_2}(n)$ up to $L$. Then for any $m_1, m_2$ and any  $L$ and $N$ sufficiently large, the following hypotheses can be considered when testing $GGC_{m_1, m_2}$ up to $N$:



\bigskip

$H_1(L, N):$	Algorithm 1/a is faster than Algorithm 1/b if and only if

\begin{equation} \label{Eqn_x}
	f_{{}_L}(m_1, m_2)<f_{{}_L}(m_2, m_1) \textrm{\; \;}  \left( \Leftrightarrow \textrm{\; \;}  
	\frac{\overline{p^*}_{{}_{L}}(m_1, m_2)}{\overline{p^*}_{{}_{L}}(m_2, m_1)}< \frac{m_2^2\varphi(m_1)}{m_1^2\varphi(m_2)}\right). 
\end{equation}

\bigskip
	
$H_2(L, N):$ Algorithm 2/a is faster than Algorithm 2/b if and only if  

\begin{equation} \label{Eqn_y}
	g_{{}_L}(m_1, m_2)<g_{{}_L}(m_2, m_1) \textrm{\; \;}   \left( \Leftrightarrow \textrm{\; \;}    \frac{\overline{p^*}_{{}_{L}}(m_1, m_2)\ln{\overline{p^*}_{{}_{L}}(m_2, m_1)}}{\overline{p^*}_{{}_{L}}(m_2, m_1)\ln{\overline{p^*}_{{}_{L}}(m_1, m_2)}} < \frac{\varphi(m_2)}{\varphi(m_1)}
	\right). 
\end{equation}

Then $H_1$ is the hypothesis that $H_1(L, N)$ is true for every $N\geq L$ where $L$ is sufficiently large. Hypothesis $H_2$ is the claim that $H_2(L, N)$ is true for every $N\geq L$ where $L$ is sufficiently large.

We tested $H_1(10^9, N_{m_1, m_2})$ and $H_2(10^9, N_{m_1, m_2})$ 
for all $489$ pairs $m_1< m_2$ relatively prime.
The pairs can be categorised as follows:

\medskip

	\textit{Group a:}  $f_{{}_{10^9}}(m_1, m_2)<f_{{}_{10^9}}(m_2, m_1)$ and $g_{{}_{10^9}}(m_1, m_2)<g_{{}_{10^9}}(m_2, m_1)$.
	
	\textit{Group b:} $f_{{}_{10^9}}(m_1, m_2)<f_{{}_{10^9}}(m_2, m_1)$ and $g_{{}_{10^9}}(m_1, m_2)>g_{{}_{10^9}}(m_2, m_1)$.

	\textit{Group c:} $f_{{}_{10^9}}(m_1, m_2)>f_{{}_{10^9}}(m_2, m_1)$ and $g_{{}_{10^9}}(m_1, m_2)<g_{{}_{10^9}}(m_2, m_1)$.

	\textit{Group d:} $f_{{}_{10^9}}(m_1, m_2)>f_{{}_{10^9}}(m_2, m_1)$ and $g_{{}_{10^9}}(m_1, m_2)>g_{{}_{10^9}}(m_2, m_1)$.
	
\medskip

Group $a$ is empty, while Groups $b, c$ and $d$ contain $227$, $258$ and $4$ pairs, respectively. For all $4$ pairs in group d at least one of the differences $|f_{10^9}(m_1, m_2)-f_{10^9}(m_2, m_1)|$ and $|g_{{}_{10^9}}(m_1, m_2)-g_{{}_{10^9}}(m_2, m_1)|$ is `small'-- less than $0,4$ -- hence it is plausibe that their group allocation may change if $L$ is sufficiently large.

Table \ref{Table_Speeds} shows the classification of the pairs  
into groups $A, B, C, D$ and $a, b, c, d$, respectively.  
In our experiment $H_1(10^9, N_{m_1, m_2})$  is true for $467$ pairs (groups $Ab, Bb, Cc$ and $Dc$),  $H_2(10^9, N_{m_1, m_2})$ holds for $476$ pairs (groups $Ac, Bb, Cc, Bd$ and $Db$) and both claims hold for $458$  pairs (groups $Bb$ and $Cc$) among all $489$ pairs. Among those $22$ pairs for which $H_1(10^9, N_{m_1, m_2})$ fails (groups $Ac, Ad, Bc, Bd$ and $Db$) in case of $15$ pairs either the running times of Algorithms 1/a and 1/b were `close' (i.e. differed by less than $10^4$ sec) or $|f_{{}_{10^9}}(m_1, m_2)-f_{{}_{10^9}}(m_2, m_1)|$ was `small' (i.e. less than $1$). For all those $13$ pairs for which $H_2(10^9, N_{m_1, m_2})$ fails (groups $Ab, Ad, Bc$ and $Dc$) either the running times of Algorithms 2/a and 2/b were 'close' (differed by less than $10^4$ sec) or $|g_{{}_{10^9}}(m_1, m_2)-g_{{}_{10^9}}(m_2, m_1)|$ was `small' (less than $1$). Hence it is plausible that for sufficiently large $N$ and $L$ the hypotheses may hold for most (or for all) of these pairs as well.



Further computational experiments and understanding the behaviours of, and developing estimations for the functions $p^*_{m_1, m_2}(n)$ could help ascertain the plausibility 
of the two hypotheses.


\subsection{Further observations regarding $p^*_{m_1, m_2}(n)$}

\begin{figure}[!htb]
	\centering
	\begin{subfigure}{0.45\textwidth}
		\includegraphics[width=\textwidth]{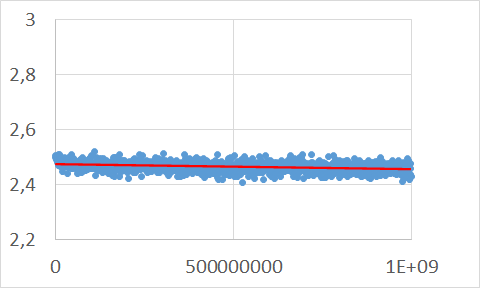}
		\caption{$m_1=1, m_2=2$}
		\label{fig:first}
	\end{subfigure}
	\hfill
	\begin{subfigure}{0.45\textwidth}
		\includegraphics[width=\textwidth]{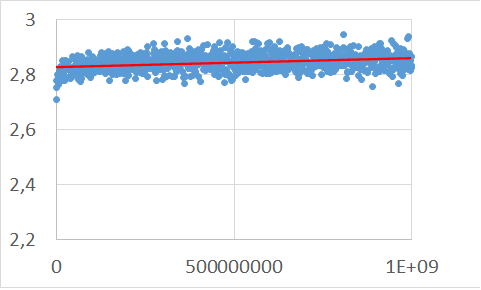}
		\caption{$m_1=2, m_2=5$}
		\label{fig:second}
	\end{subfigure}
	\hfill
	\begin{subfigure}{0.45\textwidth}
		\includegraphics[width=\textwidth]{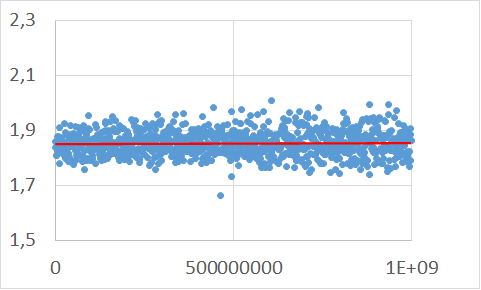}
		\caption{$m_1=23, m_2=40$}
		\label{fig:third}
	\end{subfigure}
\hfill
\begin{subfigure}{0.45\textwidth}
	\includegraphics[width=\textwidth]{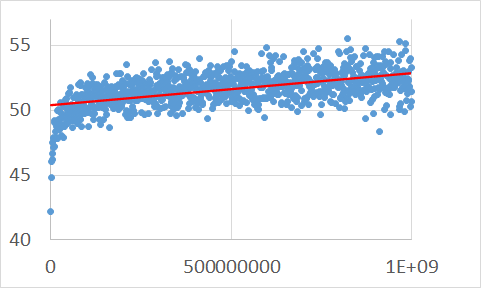}
	\caption{$m_1=1, m_2=33$}
	\label{fig:fourth}
\end{subfigure}
	
	\caption{The quotient $\frac{\textrm{average } p^*_{m_1, m_2}(n)}{\textrm{average } q^*_{m_1, m_2}(n)}$ 
		 in each interval of length $10^6$ centered at $x$, for $x=10^6k+5\cdot 10^5$ ($k=0, 1, \ldots , 10^3-1$), in cases of some  $m_1, m_2$ indicated under each subfigure.}
	\label{fig:figures}
\end{figure}


Besides their slow growths, Figure \ref{fig:Avgp} also demonstrates their closeness to smooth curves and similarity in the shapes of the graphs representing the average values of $p^*_{m_1, m_2}(n)$ in intervals of length $10^6$ up to $10^9$.

Figure \ref{fig:figures} shows the graphs of the functions $x\mapsto \textrm{average of } p^*_{m_1, m_2}(n)/\textrm{average of } q^*_{m_1, m_2}(n)$
in intervals of length $10^6$ centered at $x=10^6k+5\cdot 10^5$ ($0\leq k\leq 10^3-1$) for $m_1=1 , m_2= 2$ (Subfigure \ref{fig:first}), $m_1=2 , m_2=5 $ (Subfigure \ref{fig:second}), $m_1=23 , m_2= 40$ (Subfigure \ref{fig:third}) and $m_1=1 , m_2=33 $ (Subfigure \ref{fig:fourth}). 
The graphs -- especially the first three --  appear to be remarkably close to straight lines: the trend lines with equations $y=-2\cdot 10^{-11}x+2,4745$, $y=3\cdot 10^{-11}x+2,8297$, $y=5\cdot 10^{-12}x+1,851$ and $y=3\cdot 10^{-9}x+50,374$, indicated in Subfigures \ref{fig:first}, \ref{fig:second}, \ref{fig:third} and \ref{fig:fourth}, respectively. The values of the functions fall within the following narrow intervals between their minimum and maximum  (correct to $3$ decimal places): $[2,408 ; 2,519]$, $[2,711; 2,945]$, $[1,663; 2,01 ]$ and $ [42,200; 55,582]$ (Subfigures \ref{fig:first},\ref{fig:second}, \ref{fig:third} and \ref{fig:fourth}, respectively). If the smoothly increasing or decreasing trends of these functions continue, it suggests that the functions $L\mapsto \overline{p^*}_{{}_L}(m_1, m_2)/\overline{p^*}_{{}_L}(m_2, m_1)$
may accordingly be increasing or decreasing,  and in this case inequality \ref{Eqn_x} is either simultaneously true or false for all $L$ sufficiently large.



\nopagebreak

\section{An extension of GGC}  \label{Sect_FGGC}

Below  an extension of $GGC$  is proposed by the first author 
 -- derived in a natural way --,  which 
also generalises some other well-known conjectures regarding primes. 
It appears plausible that the requirement on $m_1, m_2$ and $n$ to be all  positive can be omitted from GGC. 
A new statement is obtained if 
$m_1$ and $m_2$ are of opposite signs (e.g. $m_1>0$ and $m_2<0$) and $n$ can be of any sign; in this case an infinite number of prime solutions in $p$ and $q$ might be possible. (Statement $(2)$ below.) 


	
	
	


\begin{EGGC}
	\label{con_FGGC}
	Let $m_1>0$ and $m_2\neq 0$ be  integers.  Then:

	\begin{enumerate}
	
\item \label{con_EGGC_1} (GGC) If $m_2>0$	then for every sufficiently large  integer $n$  satisfying the following conditions:
	
	\begin{enumerate}
		\item $\gcd(n, m_1)=\gcd(n, m_2) = \gcd(m_1, m_2)$ and   \label{con_FGGC_1a}
		\item $n\equiv m_1+m_2$   $\pmod{2^{s+1}}$, where $2^s$ is the greatest power of $2$ that is a common divisor of $m_1$ and $m_2$ \label{con_FGGC_1b} 
	\end{enumerate}
	there exist primes $p$ and $q$ such that: 
	\begin{equation} \label{Eqn_partition2}
		n=m_1p+m_2q.
	\end{equation}
\item \label{con_FGGC_2} (Extension) 	If $m_2<0$ then for every  integer $n$ satisfying Conditions \ref{con_FGGC_1a} and \ref{con_FGGC_1b},  
	equation \ref{Eqn_partition2} has infinitely many prime solutions $p$ and $q$. 
	
	\end{enumerate}
\end{EGGC}

	

 The Twin prime conjecture states that there are infinitely many twin primes, i.e. pairs of primes with a difference of $2$.
 With 
 $m_1=1$, $m_2=-1$ and $n=2$ $EGGC$ yields exactly this claim. 
  (The Twin prime conjecture is also a special case of Polignac's conjecture asserting that any positive even number $n$ can be expressed as the difference of two consecutive primes in infinitely many ways \cite{Polignac}. In its current form $EGGC$ is not a generalisation of the latter, because although with $m_1=1$, $m_2=-1$ where $n$ can be any even number it produces a similar statement, but without the condition that $p$ and $q$ are consecutive primes. A stronger version of statement (2) in $EGGC$ containing this additional requirement 
  could also be considered.)
 
 Primes of the form $2p+1$ where $p$ is also prime are called Sophie Germain primes. It is generally believed -- but has not been proved -- that there are infinitely many Sophie Germain primes \cite{Shoup}. $EGGC$ with the choice $m_1=1$, $m_2=-2$ and $n=1$ yields  exactly this assertion.

\section{Conclusions and future work} \label{Sect_Future}

The relatively small sizes of $\hat{k}_{m_1, m_2}$ in case of each pair $m_1, m_2$ tested support the plausibility of $GGC$, suggesting that the conjecture merits further investigation.

For all pairs $m_1, m_2\leq 40$ relatively prime, algorithms applying descending search were faster than those using ascending search. Heuristic arguments suggest that this is probably the case in general. However, speed rankings of the two algorithms using descending [ascending] search varied by $m_1, m_2$. The fastest algorithm can be further developed or applied potentially  in combination with sieving methods in future verification efforts. Hence it would be useful to obtain predictions for the fastest one for given $m_1, m_2$ when testing up to large tresholds. Hypotheses $H_1(10^9, N_{m_1, m_2})$ and $H_2(10^9, N_{m_1, m_2})$  were true in our implementation for most $m_1, m_2$ tested, giving support to $H_1$ and $H_2$. Further computational experiments and developing approximations to the functions $p^*_{m_1, m_2}(n)$ could help asseess their plausibility, and possibly propose  better predictions. It would  be interesting to devise  predictions for the speed rankings purely based on the values $m_1, m_2$. 

Ranking by size  of the averages $\overline{p^*}_L(m_1, m_2)$ for different $m_1, m_2\leq 40$ for $L$ sufficiently  large appears to be independent of $L$. (We could observe this in our data only when $L\leq 10^{12}$, but this is likely to be the case also for all  larger $L$.) Explaining this ranking -- and in particular, the observation that $\overline{p^*}_{10^9}(m_1, m_2)>\overline{p^*}_{10^9}(m_2, m_1)$ for every $m_1<m_2$ tested (see Table \ref{Table_AverMaxp}) -- by the properties of the numbers $m_1$ and $m_2$ is a future goal.  


Efficient sieving methods could be developed for testing $GGC$ up to higher tresholds (and potentially combined with one of the four algorithms described).




\newgeometry{top=15mm, left=10mm, right=10mm}
\section{Tables of data} \label{Sect_Tables}
\nopagebreak
\footnotesize{
\hspace*{-3cm}\begin{longtable}[c]{|c|c|c||c|c|c||c|c|c||c|c|c||c|c|c|}
	\caption{The value of $\hat{k}_{m_1, m_2}$ 
	for every  $m_1\leq m_2\leq 40$ relatively prime.} 
\label{Table_CountEx}\\
	
	\hline
	
	$\mathbf{m_1}$ & $\mathbf{m_2}$ & $\mathbf{\hat{k}_{m_1, m_2}}$ & $\mathbf{m_1}$ & $\mathbf{m_2}$ & $\mathbf{\hat{k}_{m_1, m_2}}$ & $\mathbf{m_1}$ & $\mathbf{m_2}$ & $\mathbf{\hat{k}_{m_1, m_2}}$ & $\mathbf{m_1}$ & $\mathbf{m_2}$ & $\mathbf{\hat{k}_{m_1, m_2}}$ & $\mathbf{m_1}$ & $\mathbf{m_2}$ & $\mathbf{\hat{k}_{m_1, m_2}}$  \\
	\hline
	\endhead
	
	\hline
	\endfoot
	1     & 1     & 2     & 2     & 21    & 275   & 4     & 33    & 2773  & 7     & 15    & 1192  & 9     & 25    & 6658 \\
	1     & 2     & 5     & 2     & 23    & 2209  & 4     & 35    & 9271  & 7     & 16    & 10463 & 9     & 26    & 10271 \\
	1     & 3     & 10    & 2     & 25    & 2399  & 4     & 37    & 21881 & 7     & 17    & 8104  & 9     & 28    & 6469 \\
	1     & 4     & 77    & 2     & 27    & 781   & 4     & 39    & 5443  & 7     & 18    & 6841  & 9     & 29    & 12058 \\
	1     & 5     & 24    & 2     & 29    & 4339  & 5     & 6     & 191   & 7     & 19    & 17846 & 9     & 31    & 14422 \\
	1     & 6     & 13    & 2     & 31    & 3229  & 5     & 7     & 458   & 7     & 20    & 8387  & 9     & 32    & 17021 \\
	1     & 7     & 36    & 2     & 33    & 659   & 5     & 8     & 1333  & 7     & 22    & 10729 & 9     & 34    & 14803 \\
	1     & 8     & 49    & 2     & 35    & 3733  & 5     & 9     & 274   & 7     & 23    & 13492 & 9     & 35    & 5392 \\
	1     & 9     & 28    & 2     & 37    & 11251 & 5     & 11    & 1516  & 7     & 24    & 6583  & 9     & 37    & 18976 \\
	1     & 10    & 29    & 2     & 39    & 1679  & 5     & 12    & 953   & 7     & 25    & 8618  & 9     & 38    & 21271 \\
	1     & 11    & 54    & 3     & 4     & 55    & 5     & 13    & 4582  & 7     & 26    & 22657 & 9     & 40    & 20533 \\
	1     & 12    & 25    & 3     & 5     & 62    & 5     & 14    & 3379  & 7     & 27    & 4556  & 10    & 11    & 7489 \\
	1     & 13    & 116   & 3     & 7     & 94    & 5     & 16    & 4889  & 7     & 29    & 29516 & 10    & 13    & 11051 \\
	1     & 14    & 163   & 3     & 8     & 251   & 5     & 17    & 2542  & 7     & 30    & 3217  & 10    & 17    & 13813 \\
	1     & 15    & 46    & 3     & 10    & 133   & 5     & 18    & 1187  & 7     & 31    & 25304 & 10    & 19    & 14621 \\
	1     & 16    & 473   & 3     & 11    & 140   & 5     & 19    & 3082  & 7     & 32    & 28057 & 10    & 21    & 3811 \\
	1     & 17    & 526   & 3     & 13    & 322   & 5     & 21    & 656   & 7     & 33    & 5224  & 10    & 23    & 22993 \\
	1     & 18    & 37    & 3     & 14    & 461   & 5     & 22    & 7523  & 7     & 34    & 36461 & 10    & 27    & 10537 \\
	1     & 19    & 452   & 3     & 16    & 853   & 5     & 23    & 9218  & 7     & 36    & 6091  & 10    & 29    & 28411 \\
	1     & 20    & 109   & 3     & 17    & 554   & 5     & 24    & 4229  & 7     & 37    & 39896 & 10    & 31    & 35303 \\
	1     & 21    & 88    & 3     & 19    & 616   & 5     & 26    & 16543 & 7     & 38    & 21691 & 10    & 33    & 10567 \\
	1     & 22    & 401   & 3     & 20    & 1247  & 5     & 27    & 2858  & 7     & 39    & 6472  & 10    & 37    & 45817 \\
	1     & 23    & 832   & 3     & 22    & 817   & 5     & 28    & 8237  & 7     & 40    & 30407 & 10    & 39    & 12731 \\
	1     & 24    & 97    & 3     & 23    & 2204  & 5     & 29    & 10246 & 8     & 9     & 1633  & 11    & 12    & 3623 \\
	1     & 25    & 296   & 3     & 25    & 838   & 5     & 31    & 11668 & 8     & 11    & 6509  & 11    & 13    & 13018 \\
	1     & 26    & 337   & 3     & 26    & 1777  & 5     & 32    & 12541 & 8     & 13    & 18461 & 11    & 14    & 11293 \\
	1     & 27    & 136   & 3     & 28    & 1951  & 5     & 33    & 3182  & 8     & 15    & 1399  & 11    & 15    & 1646 \\
	1     & 28    & 1157  & 3     & 29    & 1178  & 5     & 34    & 13511 & 8     & 17    & 22273 & 11    & 16    & 25723 \\
	1     & 29    & 1588  & 3     & 31    & 3358  & 5     & 36    & 4699  & 8     & 19    & 19427 & 11    & 17    & 18404 \\
	1     & 30    & 61    & 3     & 32    & 3131  & 5     & 37    & 12718 & 8     & 21    & 3517  & 11    & 18    & 6893 \\
	1     & 31    & 2918  & 3     & 34    & 1423  & 5     & 38    & 14527 & 8     & 23    & 47249 & 11    & 19    & 35254 \\
	1     & 32    & 1951  & 3     & 35    & 608   & 5     & 39    & 4954  & 8     & 25    & 14081 & 11    & 20    & 17911 \\
	1     & 33    & 214   & 3     & 37    & 3814  & 6     & 7     & 421   & 8     & 27    & 10427 & 11    & 21    & 4022 \\
	1     & 34    & 1313  & 3     & 38    & 5741  & 6     & 11    & 1361  & 8     & 29    & 43711 & 11    & 23    & 44204 \\
	1     & 35    & 226   & 3     & 40    & 2347  & 6     & 13    & 1723  & 8     & 31    & 57719 & 11    & 24    & 9707 \\
	1     & 36    & 397   & 4     & 5     & 361   & 6     & 17    & 2447  & 8     & 33    & 10841 & 11    & 25    & 31634 \\
	1     & 37    & 1616  & 4     & 7     & 1691  & 6     & 19    & 3133  & 8     & 35    & 46243 & 11    & 26    & 42073 \\
	1     & 38    & 1117  & 4     & 9     & 629   & 6     & 23    & 4901  & 8     & 37    & 57173 & 11    & 27    & 10994 \\
	1     & 39    & 272   & 4     & 11    & 2383  & 6     & 25    & 2489  & 8     & 39    & 21799 & 11    & 28    & 39167 \\
	1     & 40    & 1241  & 4     & 13    & 4073  & 6     & 29    & 10987 & 9     & 10    & 811   & 11    & 29    & 70618 \\
	2     & 3     & 17    & 4     & 15    & 1291  & 6     & 31    & 10369 & 9     & 11    & 2066  & 11    & 30    & 11021 \\
	2     & 5     & 163   & 4     & 17    & 7759  & 6     & 35    & 2059  & 9     & 13    & 3008  & 11    & 31    & 45646 \\
	2     & 7     & 89    & 4     & 19    & 12167 & 6     & 37    & 9427  & 9     & 14    & 2789  & 11    & 32    & 63601 \\
	2     & 9     & 115   & 4     & 21    & 1537  & 7     & 8     & 2711  & 9     & 16    & 7657  & 11    & 34    & 64321 \\
	2     & 11    & 673   & 4     & 23    & 24499 & 7     & 9     & 754   & 9     & 17    & 3968  & 11    & 35    & 31228 \\
	2     & 13    & 719   & 4     & 25    & 7181  & 7     & 10    & 2453  & 9     & 19    & 7498  & 11    & 36    & 18121 \\
	2     & 15    & 173   & 4     & 27    & 6511  & 7     & 11    & 2294  & 9     & 20    & 3803  & 11    & 37    & 68018 \\
	2     & 17    & 2371  & 4     & 29    & 15133 & 7     & 12    & 2371  & 9     & 22    & 11119 & 11    & 38    & 84419 \\
	2     & 19    & 1757  & 4     & 31    & 17723 & 7     & 13    & 12326 & 9     & 23    & 7454  & 11    & 39    & 26018 \\

	\newpage
	
	11    & 40    & 59399 & 15    & 22    & 8161  & 18    & 35    & 16937 & 23    & 24    & 39959 & 28    & 29    & 202273 \\
	12    & 13    & 11449 & 15    & 23    & 12428 & 18    & 37    & 53407 & 23    & 25    & 76528 & 28    & 31    & 180791 \\
	12    & 17    & 15101 & 15    & 26    & 13421 & 19    & 20    & 76319 & 23    & 26    & 106201 & 28    & 33    & 78469 \\
	12    & 19    & 8737  & 15    & 28    & 16963 & 19    & 21    & 12112 & 23    & 27    & 50872 & 28    & 37    & 250961 \\
	12    & 23    & 16739 & 15    & 29    & 29396 & 19    & 22    & 76493 & 23    & 28    & 136651 & 28    & 39    & 69259 \\
	12    & 25    & 10477 & 15    & 31    & 22636 & 19    & 23    & 110416 & 23    & 29    & 172076 & 29    & 30    & 60619 \\
	12    & 29    & 25889 & 15    & 32    & 15227 & 19    & 24    & 34129 & 23    & 30    & 26633 & 29    & 31    & 243562 \\
	12    & 31    & 18547 & 15    & 34    & 19219 & 19    & 25    & 91904 & 23    & 31    & 201812 & 29    & 32    & 370837 \\
	12    & 35    & 14303 & 15    & 37    & 21236 & 19    & 26    & 120737 & 23    & 32    & 225457 & 29    & 33    & 105254 \\
	12    & 37    & 67777 & 15    & 38    & 23873 & 19    & 27    & 26038 & 23    & 33    & 51094 & 29    & 34    & 244907 \\
	13    & 14    & 17827 & 16    & 17    & 42103 & 19    & 28    & 78671 & 23    & 34    & 163993 & 29    & 35    & 166534 \\
	13    & 15    & 3802  & 16    & 19    & 62507 & 19    & 29    & 125218 & 23    & 35    & 81274 & 29    & 36    & 97793 \\
	13    & 16    & 32507 & 16    & 21    & 12349 & 19    & 30    & 27077 & 23    & 36    & 68507 & 29    & 37    & 377122 \\
	13    & 17    & 28876 & 16    & 23    & 61861 & 19    & 31    & 169292 & 23    & 37    & 269506 & 29    & 38    & 289069 \\
	13    & 18    & 11239 & 16    & 25    & 62849 & 19    & 32    & 171469 & 23    & 38    & 273151 & 29    & 39    & 117254 \\
	13    & 19    & 30782 & 16    & 27    & 26209 & 19    & 33    & 68188 & 23    & 39    & 85906 & 29    & 40    & 228577 \\
	13    & 20    & 25913 & 16    & 29    & 133321 & 19    & 34    & 156803 & 23    & 40    & 181699 & 30    & 31    & 54337 \\
	13    & 21    & 6542  & 16    & 31    & 128783 & 19    & 35    & 69442 & 24    & 25    & 44329 & 30    & 37    & 56227 \\
	13    & 22    & 49631 & 16    & 33    & 26981 & 19    & 36    & 44647 & 24    & 29    & 83609 & 31    & 32    & 344761 \\
	13    & 23    & 44446 & 16    & 35    & 55963 & 19    & 37    & 162286 & 24    & 31    & 83507 & 31    & 33    & 87794 \\
	13    & 24    & 14221 & 16    & 37    & 186427 & 19    & 39    & 50608 & 24    & 35    & 50339 & 31    & 34    & 317567 \\
	13    & 25    & 25658 & 16    & 39    & 48067 & 19    & 40    & 103619 & 24    & 37    & 100333 & 31    & 35    & 176636 \\
	13    & 27    & 16078 & 17    & 18    & 16151 & 20    & 21    & 16129 & 25    & 26    & 110687 & 31    & 36    & 171971 \\
	13    & 28    & 74849 & 17    & 19    & 48058 & 20    & 23    & 78457 & 25    & 27    & 39586 & 31    & 37    & 363658 \\
	13    & 29    & 64634 & 17    & 20    & 37717 & 20    & 27    & 20663 & 25    & 28    & 88909 & 31    & 38    & 348349 \\
	13    & 30    & 12949 & 17    & 21    & 13382 & 20    & 29    & 142097 & 25    & 29    & 102808 & 31    & 39    & 121438 \\
	13    & 31    & 82826 & 17    & 22    & 83597 & 20    & 31    & 102659 & 25    & 31    & 165446 & 31    & 40    & 313541 \\
	13    & 32    & 80609 & 17    & 23    & 89464 & 20    & 33    & 29797 & 25    & 32    & 215743 & 32    & 33    & 108593 \\
	13    & 33    & 16024 & 17    & 24    & 39791 & 20    & 37    & 156137 & 25    & 33    & 28454 & 32    & 35    & 195197 \\
	13    & 34    & 99131 & 17    & 25    & 39332 & 20    & 39    & 26251 & 25    & 34    & 146911 & 32    & 37    & 412987 \\
	13    & 35    & 48364 & 17    & 26    & 89533 & 21    & 22    & 29191 & 25    & 36    & 87859 & 32    & 39    & 113111 \\
	13    & 36    & 31249 & 17    & 27    & 34108 & 21    & 23    & 21962 & 25    & 37    & 251206 & 33    & 34    & 136343 \\
	13    & 37    & 92006 & 17    & 28    & 51589 & 21    & 25    & 20554 & 25    & 38    & 197587 & 33    & 35    & 39994 \\
	13    & 38    & 91009 & 17    & 29    & 101834 & 21    & 26    & 33767 & 25    & 39    & 40738 & 33    & 37    & 99146 \\
	13    & 40    & 63913 & 17    & 30    & 13703 & 21    & 29    & 30746 & 26    & 27    & 39293 & 33    & 38    & 132331 \\
	14    & 15    & 2921  & 17    & 31    & 109916 & 21    & 31    & 30112 & 26    & 29    & 174451 & 33    & 40    & 71023 \\
	14    & 17    & 43423 & 17    & 32    & 120691 & 21    & 32    & 44473 & 26    & 31    & 233429 & 34    & 35    & 166597 \\
	14    & 19    & 56237 & 17    & 33    & 52004 & 21    & 34    & 47323 & 26    & 33    & 65059 & 34    & 37    & 403357 \\
	14    & 23    & 42709 & 17    & 35    & 64166 & 21    & 37    & 41794 & 26    & 35    & 142981 & 34    & 39    & 139459 \\
	14    & 25    & 23447 & 17    & 36    & 45109 & 21    & 38    & 54287 & 26    & 37    & 262897 & 35    & 36    & 52631 \\
	14    & 27    & 19787 & 17    & 37    & 203162 & 21    & 40    & 22943 & 27    & 28    & 56647 & 35    & 37    & 201062 \\
	14    & 29    & 63871 & 17    & 38    & 173681 & 22    & 23    & 108041 & 27    & 29    & 74744 & 35    & 38    & 206653 \\
	14    & 31    & 71413 & 17    & 39    & 45572 & 22    & 25    & 91277 & 27    & 31    & 54784 & 35    & 39    & 53336 \\
	14    & 33    & 19571 & 17    & 40    & 86201 & 22    & 27    & 49333 & 27    & 32    & 82343 & 36    & 37    & 113177 \\
	14    & 37    & 83717 & 18    & 19    & 35353 & 22    & 29    & 161383 & 27    & 34    & 86791 & 37    & 38    & 390367 \\
	14    & 39    & 17189 & 18    & 23    & 28153 & 22    & 31    & 133283 & 27    & 35    & 41098 & 37    & 39    & 140548 \\
	15    & 16    & 8221  & 18    & 25    & 10843 & 22    & 35    & 91579 & 27    & 37    & 94342 & 37    & 40    & 264023 \\
	15    & 17    & 6668  & 18    & 29    & 48683 & 22    & 37    & 229309 & 27    & 38    & 86143 & 38    & 39    & 188473 \\
	15    & 19    & 9664  & 18    & 31    & 37957 & 22    & 39    & 56323 & 27    & 40    & 63599 & 39    & 40    & 145279 \\

\end{longtable}

}





\newpage

{\tiny{


\begin{longtable}[c]{ccccccccccccccc}	
	
	\caption{Average and maximum values of $p^*_{m_1, m_2}(n)$ and $q^*_{m_1, m_2}(n)$ where $n\leq 10^9$ for every 
		$m_1 \leq  m_2\leq 20$ relatively prime.}
	\label{Table_AverMaxp}\\
	
	\centering

	\begin{tabular}{|l|l|l|l|l|l||l|l|l|l|l|l||l|l|l|l|l|l|} 		
		\hline
		$\mathbf{m_1}$ & 	$\mathbf{m_2}$ & \multicolumn{2}{c|}{$\mathbf{p^*_{m_1, m_2}(n)}$}  & \multicolumn{2}{c||}{$\mathbf{q^*_{m_1, m_2}(n)}$}  & 	$\mathbf{m_1}$ & $\mathbf{m_2}$ & \multicolumn{2}{c|}{$\mathbf{p^*_{m_1, m_2}(n)}$}  & \multicolumn{2}{c||}{$\mathbf{q^*_{m_1, m_2}(n)}$}  & 	$\mathbf{m_1}$ & $\mathbf{m_2}$ & \multicolumn{2}{c|}{$\mathbf{p^*_{m_1, m_2}(n)}$}  & \multicolumn{2}{c|}{$\mathbf{q^*_{m_1, m_2}(n)}$} \\
		
		 &  & \textbf{avg} & \textbf{max}  & \textbf{avg}  & \textbf{max} &   &   & \textbf{avg}  & \textbf{max}  & \textbf{avg}  & \textbf{max}  &   &   & \textbf{avg}  & \textbf{max}  & \textbf{avg}  & \textbf{max}  \\
		 
		 \hline

		1     & 1     &       &       &       &       & 4     & 9     & 241,822 & 7927  & 97,774 & 3001  & 9     & 13    & 333,584 & 10193 & 222,26 & 6761 \\
		1     & 2     & 80,839 & 3037  & 32,8  & 1609  & 4     & 11    & 494,758 & 19507 & 160,372 & 5939  & 9     & 14    & 331,513 & 10067 & 198,584 & 6337 \\
		1     & 3     & 72,911 & 2371  & 20,072 & 743   & 4     & 13    & 607,515 & 24919 & 163,502 & 6311  & 9     & 16    & 463,174 & 13627 & 241,826 & 7219 \\
		1     & 4     & 181,026 & 6971  & 32,806 & 1453  & 4     & 15    & 327,845 & 9257  & 73,338 & 2153  & 9     & 17    & 464,765 & 13007 & 227,655 & 6481 \\
		1     & 5     & 176,526 & 6833  & 26,767 & 1093  & 4     & 17    & 841,539 & 29669 & 167,531 & 6553  & 9     & 19    & 528,846 & 15649 & 229,533 & 6301 \\
		1     & 6     & 157,484 & 4969  & 17,279 & 643   & 4     & 19    & 960,026 & 32801 & 168,921 & 6947  & 9     & 20    & 449,818 & 13921 & 180,773 & 5519 \\
		1     & 7     & 281,84 & 9431  & 29,376 & 1129  & 5     & 6     & 118,745 & 3457  & 93,492 & 2801  & 10    & 11    & 370,28 & 13093 & 339,42 & 12241 \\
		1     & 8     & 393,604 & 15497 & 32,806 & 1493  & 5     & 7     & 211,95 & 8969  & 145,513 & 6871  & 10    & 13    & 454,483 & 15731 & 345,898 & 11117 \\
		1     & 9     & 245,866 & 8431  & 20,072 & 647   & 5     & 8     & 295,392 & 10369 & 172,142 & 6229  & 10    & 17    & 632,311 & 21647 & 354,305 & 14369 \\
		1     & 10    & 382,522 & 13009 & 23,958 & 1153  & 5     & 9     & 184,588 & 5333  & 97,315 & 2731  & 10    & 19    & 721,599 & 25057 & 357,248 & 13033 \\
		1     & 11    & 500,068 & 17093 & 31,678 & 1499  & 5     & 11    & 375,582 & 11839 & 156,587 & 5881  & 11    & 12    & 304,228 & 8821  & 267,184 & 8293 \\
		1     & 12    & 342,648 & 11261 & 17,279 & 673   & 5     & 12    & 257,838 & 7309  & 93,511 & 2969  & 11    & 13    & 542,8 & 17299 & 452,035 & 16829 \\
		1     & 13    & 612,063 & 23663 & 32,294 & 1297  & 5     & 13    & 459,273 & 16477 & 159,551 & 5521  & 11    & 14    & 542,423 & 20359 & 406,549 & 15227 \\
		1     & 14    & 611,042 & 20359 & 26,557 & 1129  & 5     & 14    & 459,822 & 15773 & 141,315 & 4651  & 11    & 15    & 295,49 & 8941  & 203,796 & 5527 \\
		1     & 15    & 332,373 & 9127  & 15,379 & 557   & 5     & 16    & 638,409 & 24677 & 172,167 & 6451  & 11    & 16    & 754,067 & 26839 & 494,633 & 17863 \\
		1     & 16    & 849,623 & 33997 & 32,803 & 1597  & 5     & 17    & 637,215 & 22751 & 163,446 & 5657  & 11    & 17    & 751,415 & 25621 & 463,029 & 17713 \\
		1     & 17    & 846,422 & 32779 & 33,084 & 1381  & 5     & 18    & 398,036 & 10499 & 93,511 & 2963  & 11    & 18    & 470,391 & 14251 & 267,222 & 7681 \\
		1     & 18    & 529,975 & 15313 & 17,28 & 701   & 5     & 19    & 726,834 & 25609 & 164,784 & 5711  & 11    & 19    & 856,789 & 27581 & 466,872 & 15467 \\
		1     & 19    & 964,977 & 33791 & 33,364 & 1321  & 6     & 7     & 149,317 & 4597  & 129,94 & 3923  & 11    & 20    & 734,055 & 26497 & 370,239 & 12853 \\
		1     & 20    & 825,834 & 29209 & 23,957 & 1069  & 6     & 11    & 267,139 & 8543  & 139,856 & 4813  & 12    & 13    & 328,85 & 9871  & 310,206 & 9479 \\
		2     & 3     & 69,352 & 2083  & 43,626 & 1399  & 6     & 13    & 328,759 & 10883 & 142,486 & 4957  & 12    & 17    & 459,61 & 13033 & 317,598 & 10657 \\
		2     & 5     & 172,137 & 6379  & 60,482 & 2459  & 6     & 17    & 459,546 & 14731 & 145,948 & 4201  & 12    & 19    & 523,692 & 14699 & 320,198 & 9437 \\
		2     & 7     & 277,107 & 12011 & 66,282 & 2663  & 6     & 19    & 523,593 & 16703 & 147,126 & 4423  & 13    & 14    & 552,943 & 19889 & 499,815 & 16843 \\
		2     & 9     & 241,78 & 7129  & 43,628 & 1549  & 7     & 8     & 323,92 & 12589 & 277,119 & 11197 & 13    & 15    & 301,049 & 8539  & 250,574 & 7151 \\
		2     & 11    & 494,633 & 21107 & 71,487 & 3061  & 7     & 9     & 202,501 & 5717  & 154,108 & 4271  & 13    & 16    & 768,659 & 28463 & 607,268 & 25127 \\
		2     & 13    & 607,339 & 21383 & 72,924 & 3049  & 7     & 10    & 315,347 & 9769  & 207,433 & 6841  & 13    & 17    & 765,986 & 25747 & 566,894 & 21851 \\
		2     & 15    & 327,714 & 9049  & 32,917 & 1031  & 7     & 11    & 411,838 & 15131 & 249,973 & 9439  & 13    & 18    & 479,435 & 15199 & 328,849 & 9277 \\
		2     & 17    & 841,438 & 30859 & 74,71 & 3121  & 7     & 12    & 282,761 & 9137  & 149,358 & 4663  & 13    & 19    & 873,509 & 33703 & 571,528 & 22079 \\
		2     & 19    & 959,87 & 34039 & 75,341 & 3001  & 7     & 13    & 504,267 & 18593 & 254,754 & 10099 & 13    & 20    & 748,115 & 27953 & 454,483 & 14851 \\
		3     & 4     & 97,757 & 2939  & 69,363 & 2411  & 7     & 15    & 274,865 & 7499  & 116,406 & 3583  & 14    & 15    & 270,331 & 7789  & 248,994 & 6689 \\
		3     & 5     & 97,292 & 2909  & 55,338 & 1709  & 7     & 16    & 700,594 & 22783 & 277,12 & 10357 & 14    & 17    & 693,436 & 25121 & 566,271 & 20717 \\
		3     & 7     & 154,073 & 4517  & 60,42 & 1789  & 7     & 17    & 698,137 & 24109 & 260,916 & 11069 & 14    & 19    & 791,453 & 27277 & 570,866 & 20873 \\
		3     & 8     & 211,872 & 6869  & 69,37 & 2383  & 7     & 18    & 436,631 & 13367 & 149,359 & 4481  & 15    & 16    & 347,11 & 9521  & 327,84 & 8893 \\
		3     & 10    & 208,887 & 6359  & 51,951 & 1471  & 7     & 19    & 796,433 & 27583 & 263,145 & 10289 & 15    & 17    & 349,899 & 9539  & 308,256 & 8179 \\
		3     & 11    & 271,626 & 8231  & 64,839 & 2113  & 7     & 20    & 682,396 & 23689 & 207,435 & 6841  & 15    & 19    & 398,729 & 10979 & 310,687 & 9109 \\
		3     & 13    & 333,472 & 10733 & 66,064 & 1999  & 8     & 9     & 241,796 & 7027  & 211,92 & 6961  & 16    & 17    & 841,42 & 30727 & 787,381 & 29531 \\
		3     & 14    & 331,38 & 10259 & 57,045 & 1867  & 8     & 11    & 494,594 & 18481 & 348,832 & 13499 & 16    & 19    & 959,865 & 35327 & 793,796 & 28631 \\
		3     & 16    & 463,076 & 13553 & 69,361 & 2239  & 8     & 13    & 607,287 & 23887 & 355,623 & 12107 & 17    & 18    & 491,044 & 14149 & 459,684 & 12953 \\
		3     & 17    & 464,638 & 12503 & 67,628 & 2269  & 8     & 15    & 327,799 & 9091  & 158,333 & 4817  & 17    & 19    & 894,547 & 33721 & 790,894 & 26927 \\
		3     & 19    & 528,697 & 15217 & 68,167 & 2063  & 8     & 17    & 841,438 & 31081 & 364,403 & 15749 & 17    & 20    & 765,917 & 28429 & 632,339 & 25237 \\
		3     & 20    & 449,579 & 12659 & 51,956 & 1579  & 8     & 19    & 959,894 & 42727 & 367,416 & 13999 & 18    & 19    & 523,747 & 14897 & 495,035 & 16943 \\
		4     & 5     & 172,187 & 7109  & 135,388 & 5521  & 9     & 10    & 208,976 & 6469  & 180,762 & 5501  & 19    & 20    & 772,014 & 28729 & 721,715 & 24071 \\
		4     & 7     & 277,169 & 11497 & 148,746 & 5939  & 9     & 11    & 271,709 & 8363  & 218,093 & 6827  &       &       &       &       &       &  \\
		\hline
	\end{tabular}%
	\label{tab:addlabel}%
\end{longtable}%
}
}
	
\nopagebreak

\begin{table}[htbp]
	\centering
	\caption{The five greatest, smallest and the average values of $\max{p^*_{m_1, m_2}}$ and of $\overline{p^*}_{m_1, m_2}$ up to $10^9$ and of $\hat{k}_{m_1, m_2}$ over all pairs $m_1, m_2\leq 40$ relatively prime. 
	}
	\begin{tabular}{|l|r|c|c|c|c|}
		\hline
		& \multicolumn{2}{c|}{\textbf{5 greatest values}}  &	    \multicolumn{2}{c|}{\textbf{5 smallest values}} &   \textbf{Average} \\

		&  \multicolumn{1}{c|}{\textbf{value}}     & $\mathbf{(m_1, m_2)}$ &  \multicolumn{1}{c|}{\textbf{value}} & $\mathbf{(m_1, m_2)}$ & \textbf{value}\\
		\hline
		\hline	
		
		$\mathbf{max \; {p^*_{m_1, m_2}}}$ \textbf{up to} $\mathbf{10^9}$ & $78697$ & $(32, 37)$ &  $449$ & $(30, 1)$ &  $22889,33538$\\
		& $77723$ & $(23, 37)$  & $557$ & $(17, 1)$  & \\
		
		& $77267$ & $(37, 38)$ & $571$ & $(39, 1)$ & \\
		
		& $76379$ & $(29, 38)$ & $599$ & $(21, 1)$ & \\
		
		& $75989$ & $(1, 38)$  & $631$ & $(24, 1)$  & \\
		\hline		
		
		$\mathbf{\overline{p^*}_{{}_{m_1, m_2}}}$ \textbf{up to} $\mathbf{10^9}$ & 2064,47552 & $(1, 37)$ & $12,74269$ & $(30, 1)$  & $687,7063317$\\
		
		& $2059,89836$ & $(1, 38)$ & $15,37864$ & $(15, 1)$  & \\
		
		& $2059,17801$ & $(16, 37)$  & $16,68819$ & $(21, 1)$  & \\
		
		& $2059,1531$ & $(32, 37)$  & $17,27778$ & $(36, 1)$  & \\
		
		& $2058,97664$ & $(2, 37)$  & $17,27898$ & $(6, 1)$  & \\
		
		\hline

		$\mathbf{\hat{k}_{m_1, m_2}}$ 	& $412987$ & $(32, 37)$, $(37, 32)$ & $2$ & $(1, 1)$  & $52004,838776$\\
		

		& $403357$ & $(34, 37)$, $(37, 34)$  & $5$ & $(1, 2)$, $(2, 1)$  & \\
		

		& $390367$ & $(37, 38)$, $(38, 37)$ & $10$ & $(1, 3)$, $(3, 1)$  & \\
		

		& $377122$ & $(29, 37)$,  $(37, 29)$ & $13$ & $(1, 6)$, $(6, 1)$  & \\
		
		
		& $370837$ & $(29, 32)$, $(32, 29)$  & $17$ & $(2, 3)$, $(3, 2)$  & \\
		
		
		\hline	
	\end{tabular}%
	\label{tab:addlabel}%
	
\end{table}%

\bigskip

{\footnotesize{
		\begin{table}
			
			\caption{Classification of all pairs $ m_1<m_2\leq 40$ relatively prime into groups $A, B, C, D$ and $a, b, c, d$, 
				indicated in the first column by upper and lower case letters, respectively. 
			} \label{Table_Class}
			\medskip
			
			\begin{tabularx}{\textwidth}{|p{1cm}|X|}
				
				\hline
				\textbf{Group} &  \textbf{The ordered pairs } $\mathbf{(m_1, m_2)}$ \textbf{ contained by the group} \\	
				\hline
				\hline	
				
				$\mathbf{Ab:}$ & (1, 3)
				(1, 9),
				(1, 15), (1, 21), (1, 33), (1, 39)
				\\
				\hline
				
				$\mathbf{Ac:}$ & (1, 7),
				(1, 11),
				 (1, 13), (1, 17), (1, 19), (1, 25), (1, 31), (1, 37), (2, 9), (2, 15), (2, 21), (7, 11) \\
				\hline
				
				$\mathbf{Ad:}$ & (1, 5), 
				(1, 27), (1, 35) \\
				\hline
				
				$\mathbf{Bb:}$ & (1, 2), (1, 4), (1, 6), (1, 8), (1, 10), (1, 12), (1, 14), (1, 16), (1, 18), (1, 20), (1, 22), (1, 24), (1,~ 26), (1, 28), (1, 30), (1, 32), (1, 34), (1, 36), (1, 38), (1, 40), (3, 4), (3, 8), (3, 10), (3, 14), (3, 16), (3, 20), (3, 22), (3, 26), (3, 28), (3, 34), (3, 38), (3, 40), (5, 6), (5, 8), (5, 9), (5, 12), (5, ~ 14), (5, 16), (5, 18), (5, 21), (5, 22), (5, 24), (5, 26), (5, 27), (5, 28), (5, 32), (5, 33), (5,~ 34), (5, 36), (5, 38), (5, 39), (7, 8), (7, 9), (7, 10), (7, 12), (7, 15), (7, 16), (7, 18), (7, 20), (7, 22), (7, 24), (7, 26), (7, 27), (7, 30), (7, 32), (7, 33), (7, 34), (7, 36), (7, 38), (7, 39), (7, 40), (9,~ 10), (9, 14), (9, 16), (9, 20), (9, 22), (9, 26), (9, 28), (9, 34), (9, 38), (9, 40), (11, 12), (11, 14), (11,~ 15), (11, 16), (11, 18), (11, 20), (11, 21), (11, 24), (11, 25), (11, 26), (11, 27), (11, 28), (11, ~ 30), (11, 32), (11, 34), (11, 35), (11, 36), (11, 38), (11, 39), (11, 40), (13, 14), (13, 15), (13, 16), (13, 18), (13, 20), (13, 21), (13, 22), (13, 24), (13, 25), (13, 27), (13, 28), (13, 30), (13, 32), (13, 33), (13, 34), (13, 35), (13, 36), (13, 38), (13, 40), (15, 16), (15, 22), (15, 26), (15, 28), (15, 34), (15, 38), (17, 18), (17, 20), (17, 21), (17, 22), (17, 24), (17, 25), (17, 26), (17, 27), (17, 28), (17, 30), (17, 32), (17, 33), (17, 35), (17, 36), (17, 38), (17, 39), (17, 40), (19, 20), (19, 21), (19, 22), (19, 24), (19, 25), (19, 26), (19, 27), (19, 28), (19, 30), (19, 32), (19, 33), (19, 34), (19, 35), (19, 36), (19, 39), (19, 40), (21, 22), (21, 26), (21, 32), (21, 34), (21, 38), (21, 40), (23, 24), (23, 25), (23, 26), (23, 27), (23, 28), (23, 30), (23, 32), (23, 33), (23, 34), (23, 35), (23, 36), (23, 38), (23, 39), (25, 26), (25, 27), (25, 28), (25, 32), (25, 33), (25, 34), (25, 36), (25, 38), (25, 39), (27, 28), (27, 32), (27, 34), (27, 38), (27, 40), (29, 30), (29, 32), (29, 33), (29, 34), (29, 35), (29, 36), (29, 38), (29, 39), (31, 32), (31, 33), (31, 34), (31, 35), (31, 36), (31, 38), (31, 39), (33, 34), (33, 38), (33, 40), (35, 36), (35, 38), (35, 39), (37, 38), (37, 39), (39, 40) \\
				\hline
				
					$\mathbf{Bc:}$ & (15, 32) 
				\\
				\hline
				
				$\mathbf{Bd:}$ & (3, 32)  \\
				\hline

				$\mathbf{Cc:}$ & (1, 23), (1, 29), (2, 3), (2, 5), (2, 7), (2, 11), (2, 13), (2, 17), (2, 19), (2, 23), (2, 25), (2, 27), (2, 29), (2, 31), (2, 33), (2, 35), (2, 37), (2, 39), (3, 5), (3, 7), (3, 11), (3, 13), (3, 17), (3, 19), (3, 23), (3, 25), (3, 29), (3, 31), (3, 35), (3, 37), (4, 5), (4, 7), (4, 9), (4, 11), (4, 13), (4, 15), (4, 17), (4, 19), (4, 21), (4, 23), (4, 25), (4, 27), (4, 29), (4, 31), (4, 33), (4, 35), (4, 37), (4, 39), (5, 7), (5, 11), (5, 13), (5, 17), (5, 19), (5, 23), (5, 29), (5, 31), (5, 37), (6, 7), (6, 11), (6, 13), (6, 17), (6, 19), (6, 23), (6, 25), (6, 29), (6, 31), (6, 35), (6, 37), (7, 13), (7, 17), (7, 19), (7, 23), (7, 25), (7, 29), (7, 31), (7, 37), (8, 9), (8, 11), (8, 13), (8, 15), (8, 17), (8, 19), (8, 21), (8, 23), (8, 25), (8, 27), (8, 29), (8, 31), (8, 33), (8, 35), (8, 37), (8, 39), (9, 11), (9, 13), (9, 17), (9, 19), (9, 23), (9, 25), (9, 29), (9, 31), (9, 35), (9, 37), (10, 11), (10, 13), (10, 17), (10, 19), (10, 21), (10, 23), (10, 27), (10, 29), (10, 31), (10, 33), (10, 37), (10, 39), (11, 13), (11, 17), (11, 19), (11, 23), (11, 31), (11, 37), (12, 13), (12, 17), (12, 19), (12, 23), (12, 25), (12, 29), (12, 31), (12, 35), (12, 37), (13, 17), (13, 19), (13, 23), (13, 29), (13, 31), (13, 37), (14, 15), (14, 17), (14, 19), (14, 23), (14, 25), (14, 27), (14, 29), (14, 31), (14, 33), (14, 37), (14, 39), (15, 17), (15, 19), (15, 23), (15, 29), (15, 31), (15, 37), (16, 17), (16, 19), (16, 21), (16, 23), (16, 25), (16, 27), (16, 29), (16, 31), (16, 33), (16, 35), (16, 37), (16, 39), (17, 23), (17, 29), (17, 31), (17, 37), (18, 19), (18, 23), (18, 25), (18, 29), (18, 31), (18, 35), (18, 37), (19, 23), (19, 29), (19, 31), (19, 37), (20, 21), (20, 23), (20, 27), (20, 29), (20, 31), (20, 33), (20, 37), (20, 39), (21, 23), (21, 25), (21, 29), (21, 31), (21, 37), (22, 23), (22, 25), (22, 27), (22, 29), (22, 31), (22, 35), (22, 37), (22, 39), (23, 31), (23, 37), (24, 25), (24, 29), (24, 31), (24, 35), (24, 37), (25, 29), (25, 31), (25, 37), (26, 27), (26, 29), (26, 31), (26, 33), (26, 35), (26, 37), (27, 29), (27, 31), (27, 35), (27, 37), (28, 29), (28, 31), (28, 33), (28, 37), (28, 39), (29, 31), (29, 37), (30, 31), (30, 37), (31, 37), (32, 33), (32, 35), (32, 37), (32, 39), (33, 35), (33, 37), (34, 35), (34, 37), (34, 39), (35, 37), (36, 37), (38, 39) \\
				\hline

				$\mathbf{Db:}$ & (9, 32),
				(23, 40), (29, 40), (31, 40), (37, 40)  \\
				\hline

				\textbf{Dc:} & (11, 29), 
				(17, 19), 
				 (23, 29)
				\\
				\hline
			\end{tabularx}
			\label{Table_Speeds} \end{table}
}}

\clearpage

\appendix

\section{Pseudocodes} \label{Sect_Pseudo}

{\begin{small}
		\begin{algorithm}[H] \label{Funct_Genm1pr}
			\Fn{\Gmpr{$\alpha$}}{
				\Input{array $isprime$}
				
				\tcc{\textbf{Global variables used: $m_1, m_2$ and $\alpha$}}
				\Output{for every $0\leq r<m_2$ 
					array $\mathit{m_1p}[r]$  containing all numbers of the form $m_1p$ where $p$ is prime such that $m_1p\leq \alpha$ and $r= m_1p \mod{m_2}$}
					
					
					$inc \leftarrow 2m_1\mod{m_2}$\;
					
					add($m_1p[inc]$, $2m_1$)\;
					$L\leftarrow length(isprime)$\;
					
					$r \leftarrow 3m_1\mod{m_2}$\;
					$j\leftarrow 0$\;

					\While{$j\leq L-1$ and $m_1(2j+3)\leq \alpha$ }{
						
						\If{$isprime[j]=1$}{

							add($m_1p[r]$, $m_1\cdot (2j+3)$)\;			
							
						}
						$r \leftarrow (r \geq m_2-inc)$ ? $r+inc-m_2$ : $r+inc$\;
						
						$j \leftarrow j+1$\;
					}

				}
		\end{algorithm}
\end{small}}

\bigskip

{\begin{small}
		\begin{algorithm}[H]
			\Fn{\Gmqr{$C,D$}}{
				\Input{integers $0 \leq C<D$  such that $2m_1m_2| C$ and $2m_1m_2 |D$.}

				\tcc{Global variables used: $m_1, m_2$ and array $primes$}

				\Output{arrays \textbf{m\textsubscript{2}q[$r$]} ($0\leq r<m_1$) containing all numbers of the form $m_2q$ in interval $[C, D)$ where $q$ is prime and  $r = m_2q \mod{m_1}$.}
				\tcc{Initialization}
				$c \leftarrow \frac{C}{m_2}; d \leftarrow \frac{D}{m_2}; u \leftarrow \frac{d-c}{2}-1$ \;
				b$[i] \leftarrow 1, (0\leq i \leq u, \mbox{ if } c=0 \mbox{ then  } b[0] \leftarrow 0)$\; 
				$j \leftarrow 1$ \;
				\tcc{Sieving odd numbers in $[c, d)$ using odd primes, in order to produce boolean array $b$ such that for every $0\leq i\leq u$: $b[i]=1$ if and only if $c+2i+1$ is prime.}

				\While{$p=\mbox{primes}[j]<\sqrt{d}$}{ 
					
					\tcc{Setting starting point for sieving with first/next prime $p$}
					\eIf{$p \geq \sqrt{c}$}{$s \leftarrow p^2$\;}
					{$s \leftarrow 2p\cdot \left( \lfloor \frac{c+p-1}{2p}\rfloor\right)+p $\;}
					$k \leftarrow \frac{s-c-1}{2}$\;
					
					\tcc{Sieving with prime $p$}
					\While{$k\leq u$}{
						$b[k] \leftarrow 0$\;
						$k \leftarrow k+p$\;
					}
					$j \leftarrow j+1$\;			
				}
				
				\tcc{Populating arrays $m_2q[r]$}

				\tcc{Handling special case when $c\leq 2<d$}
				
				\If{$c\leq 2$ and $d>2$}{
					
					$r \leftarrow 2m_2  \mod{m_1}$\;
					
					add$(m_2q[r], 2m_2)$\;					
				}

				\tcc{Populating arrays $m_2q[r]$  ($0\leq r<m_1$) with values $m_2q$ where $c\leq q <d$ is odd prime such that $r = m_2q \mod{m_1}$}

				$i \leftarrow 0$\;
				$r \leftarrow m_2  \mod{m_1}$\;
				$inc \leftarrow 2m_2  \mod{m_1}$\;
				\While{$i \leq u$}{
					\If{$b[i]=1$}{add($m_2q[r]$, $m_2(c+2i+1)$)\;}
					$i \leftarrow i+1$\;
					$r \leftarrow (r \geq m_1-inc)$ ? $r+inc-m_1$ : $r+inc$\;
				}			
				
			}
		\end{algorithm}
\end{small}}

\bigskip

{\begin{small}
		\begin{algorithm}[H]
			\Fn{\Gismq{$C,D$}}{
				\Input{integers $0\leq C<D$ such that $2m_2| C$ and $2m_2 |D$.}
				
				\tcc{Global variables used: $m_1, m_2$ and array $primes$}

				\Output{boolean array $ism_2q$ of length $D-C$ such that for every $0\leq i\leq D-C-1$: $ism_2q[i]=1$ if and only if $C+i=m_2q$ for some prime $q$. 
				}
				\tcc{Initialization}
				$c \leftarrow \frac{C}{m_2}; d \leftarrow \frac{D}{m_2}; u \leftarrow \frac{d-c}{2}-1$\;
				$b[i] \leftarrow 1, (0\leq i \leq u, \mbox{ if } c=0 \mbox{ then  } b[0] \leftarrow 0)$\; 
				$j \leftarrow 1$\;
				\tcc{Sieving odd numbers in $[c, d)$ using odd primes, in order to produce boolean array $b$ such that for every $0\leq i\leq u$: $b[i]=1$ if and only if $c+2i+1$ is prime.}

				\While{$p:=\mbox{primes}[j]<\sqrt{d}$}{ 
					\tcc{Setting starting point for sieving with new prime $p$}
					\eIf{$p \geq \sqrt{c}$}{$s \leftarrow p^2$\;}
					{$s \leftarrow 2p\cdot \left( \lfloor \frac{c+p-1}{2p}\rfloor\right)+p $\;}
					$k \leftarrow \frac{s-c-1}{2}$\;
					
					\tcc{Sieving with prime $p$}
					
					\While{$k\leq u$}{
						$b[k] \leftarrow 0$\;
						$k \leftarrow k+p$\;
					}
					$j \leftarrow j+1$\;			
				}
				
				\tcc{Preparing array $ism_2q$}
				
				\tcc{Initialization}
				
				\For{$i=0$ to  $D-C-1$}{
					
					$ism\textsubscript{2}q[i] \leftarrow 0$\;
					
				}

				\tcc{Handling special case when $c\leq 2<d$}

				\If{$c\leq 2$ and $d>2$}{
					
					$ism_2q[2m_2-C] \leftarrow 1$\;
					
				}
				
				\tcc{Setting values in $ism_2q$}
				\For{$i  =0$ to $u$}{
					\If{$b[i]=1$}{$ism\textsubscript{2}q[2m_2i+m_2] \leftarrow 1$\;}
				}			
				
			}
		\end{algorithm}
\end{small}}

\bigskip

{\begin{small}
		\begin{algorithm}[H]
			\Fn{\Checkone{$A,B$}}{
				\Input{
					\begin{itemize}
						
						\item	 integers $0\leq A<B$ such that $2m_1m_2|A$ and $2m_1m_2|B$ and 
						
						\item arrays $m_2q[r]$ for every $0\leq r<m_1$, containing all numbers of the form  $m_2q$ in the \\ \noindent interval  $[\max{\{0, A-\alpha\}}, B)$ where $q$ is prime such that $r= m_2q \mod m_1$.
					\end{itemize}
					
					\tcc{Global variables used: $m_1, m_2$, arrays $res$ and $ism_1p$.}
					
					\Output{array $residual$ containing all numbers $n$ in interval $[A, B)$ satisfying the conditions of $GGC_{m_1, m_2}$ for which there do not exist primes $p$ and $q$ such that $n=m_1p+m_2q$ and $m_1p\leq \alpha$.}
					\tcc{Initialization}
					

						$r \leftarrow (m_1+m_2 \textrm{ is even})$ ? $m_1: m_1+1$; $s \leftarrow (m_1+m_2 \textrm{ is even})$ ? $lcm_{m_1, m_2}:lcm_{m_1, m_2}+1$; $n \leftarrow (m_1+m_2 \textrm{ is even})$ ? $B: B+1$\;
					
					$l[j] \leftarrow \len(m\textsubscript{2}q[j]), (0\leq j < m_1)$\; 

					\tcc{\textbf{Process}}
					\While{$n>A+1$}{

						$n \leftarrow n-2$\;
						$r \leftarrow (r<2)$ ? $r+m_1-2$ : $r-2$\;
						$s \leftarrow (s<2)$ ? $s+lcm_{m_1, m_2}-2$ : $s-2$\;
						
						\tcc{If $n$ satisfies conditions of $GGC_{m_1, m_2}$ search for $m_2q^{**}_{m_1, m_2}(n)$.}
						\If{$res[s]=1$}{
							
							\tcc{Setting starting point in $m_2q[r]$ for search for $m_2q^{**}_{m_1, m_2}(n)$.}
							\While{$l[r]>0$ and $n<m_2q[r][l[r]-1]+2m_1$}{					
								$l[r]\leftarrow l[r]-1$\;					
							}

							\eIf{$l[r]=0$}{
								add($residual$, $n$)\;
							}{
								
								$i\leftarrow l[r]-1$\;


								\tcc{Search for $m_2q^{**}_{m_1, m_2}(n)$ starts.}	
								
								\While{$i\geq 0$}{
									
									\tcc{$m_2q^{**}_{m_1, m_2}(n)$ has not been found and $m_2q[r][i]$ has become too small.}
									
									\uIf{$n-m_2q[r][i]\geq length(ism_1p)$}{
										
										add($residual$, $n$)\;
										
										break\;
									}

									\tcc{$m_2q[r][i]=m_2q^{**}_{m_1, m_2}(n)$; $p^*_{m_1, m_2}(n)$ and $q^{**}_{m_1, m_2}(n)$ optionally can be saved.}
									
									\uElseIf{$ism_1p[n-m_2q[r][i]]=1$}{
										$p^*_{m_1, m_2}(n)\leftarrow \frac{n-m_2q[r][i]}{m_1}$\;
										
										$q_{m_1,m_2}^{**}(n)\leftarrow \frac{m_2q[r][i]}{m_2}$\;
										
										
										break\;								
									}
									
									\tcc{$m_2q[r][i] \neq m_2q^{**}_{m_1, m_2}(n)$ and $m_2q[r][i]$ is not too small yet.}
									
									\uElse{$i\leftarrow i-1$\;}
								}
								
								\If{$i=-1$}{
									
									add($residual$, $n$)\;	}

							}				
						}			
					}
				}
			} 
		\end{algorithm}
		
\end{small}}

\bigskip

\begin{algorithm}[H]
	\Fn{\Checktwo{$A,B$}}{
		\Input{
			\begin{itemize}	
				
				\item	integers $0\leq A<B$  such that $2m_1m_2|A$ and $2m_1m_2|B$ and
				
				\item boolean array $ism_2q$ of length $B-C$ , where $C=\max{\{0, A-\alpha\}}$, such that \\ \noindent for every $0\leq i<B-C$: $ism_2q[i]=1$ if and only if $C+i=m_2q$ for some prime $q$.
				
			\end{itemize}	
		}
		
		\tcc{Global variables used: $m_1, m_2$ and arrays $res$ and $m_1p[r]$ for every $0\leq r<m_2$.}
		\Output{array \textit{residual} containing all numbers $n$ in interval $[A, B)$ satisfying the conditions of $GGC_{m_1, m_2}$ for which there do not exist primes $p$ and $q$ such that $n=m_1p+m_2q$ and $m_1p\leq \alpha$.}
		\tcc{Initialization}

		$r \leftarrow (m_1+m_2 \textrm{ is even})$ ? $m_2: m_2+1$; $s \leftarrow (m_1+m_2 \textrm{ is even})$ ? $lcm_{m_1, m_2}:lcm_{m_1, m_2}+1$; $n \leftarrow (m_1+m_2 \textrm{ is even})$ ? $B: B+1$\;

		$l[j] \leftarrow \len(m\textsubscript{1}p[j])$ for every $0\leq j < m_2$\; 
		\tcc{\textbf{Process}}
		\While{$n>A+1$}{

			$n \leftarrow n-2$\;
			$r \leftarrow (r<2)$ ? $r+m_2-2$ : $r-2$\;
			$s \leftarrow (s<2)$ ? $s+lcm_{m_1, m_2}-2$ : $s-2$\;

			\tcc{If $n$ satisfies conditions of $GGC_{m_1, m_2}$ search for $m_1p^{*}_{m_1, m_2}(n)$.}

			\If{$res[s]=1$}{
				\For{$i=0$ to $l[r]-1$}{			
					\If{$n-m_1p[r][i]\geq C$ and $ism_2q\left[ n-m_1p[r][i]-C\right]=1 $}{
						$p^*_{m_1,m_2}(n) \leftarrow \frac{m_1p[r][i]}{m_1}$\;
						$q^{**}_{m_1,m_2}(n) \leftarrow \frac{n-m_1p[r][i]}{m_2}$\;
						break;
					}
					\If{$i=l[r]-1$ or $n-m_1p[r][i]< C$}{
						add($residual$, $n$)\;
						break; }	
					
				}
			}			
		}
	}
\end{algorithm}

{\begin{small}
		\begin{algorithm}[H] \label{Fn_GGCone}
			\Fn{\GGCone{$N,m_1,m_2, \triangle$, $\alpha$}}{
				\Input{positive integers $N, m_1,m_2$, $\triangle$ and $\alpha$  such that $\gcd (m_1, m_2)=1$, $N > 9$, $2m_1m_2|N$, $2m_1m_2| \triangle$ and  
					$\alpha\leq \triangle$.
				}
				\Output{array \textit{residual} containing all numbers $n\leq N$ satisfying the conditions  of $GGC_{m_1, m_2}$ for which there do not exist primes $p$ and $q$ such that $n=m_1p+m_2q$ and $m_1p\leq \alpha$.}
				
				\tcc{\textbf{Start Phase I: Unsegmented phase}}
				\tcc{Generating array 
					$primes$.} 
				\textbf{SmallPrimes($\max{(\lfloor \sqrt{\frac{N_{m_1, m_2}}{m_2}} \rfloor, \lfloor\frac{\alpha}{m_1}\rfloor )})$}\;

				\tcc{Assigning values to array  $ism_1p$.}

				\textbf{GenerateIsm1p}($\alpha$)\;
				
				\tcc{Assigning values to array  $res$.}
				
				\textbf{GenerateResiduePattern}($m_1, m_2$)\;

				\tcc{\textbf{Start Phase II: Segmented phase}}
				
				\tcc{\textbf{Initialization}}

				Set arrays $residual$ and $m_2q[r]$ ($0\leq r< m_1$) empty\;
				
				$A\leftarrow 0$\;

				\tcc{\textbf{Start segmented computation}}
				
				\While{$A<N$}{
					
					$B\leftarrow \min{\{A+\triangle, N\}}$\;

					\tcc{Keeping only those values in each array $m_2q[r]$ generated in previous iteration which are greater than $A-\alpha$ and removing all other values.}	
					
					\If{$A>0$}{
						\For{$r=0$ \KwTo $m_1-1$}{

							$i\leftarrow 0$\;
							
							\While{$i<length(m_2q[r])$ and $m_2q[r][i]< A-\alpha$}{		        
								$i\leftarrow i + 1$\;
							}
							
							\If{$i\neq 0$}{
								remove\_interval($mq_2[r]$, $[0\ldots i-1]$)}
						}		    
					}
					
					\tcc{Assigning new values to arrays $m_2q[r]$.}
					\textbf{Generatem2qr$(A,B)$}\;
					\tcc{Checking $GGC_{{}_{m_1,m_2}}$ in new interval.}
					
					\textbf{Check1$(A,B)$}\;
					
					$A\leftarrow A+\triangle$\;

				}		
			}
		\end{algorithm}
		
\end{small}}




{\begin{small}
		\begin{algorithm}[H] \label{Fn_GGCtwo}
			\Fn{\GGCtwo{$N,m_1,m_2, \triangle$, $\alpha$}}{
				\Input{$N, m_1,m_2$, $\triangle$ and $\alpha$ positive integers such that $\gcd (m_1, m_2)=1$, $N >9$, $2m_1m_2|N$, $2m_1m_2| \triangle$ and 
					$\alpha\leq \triangle$.}
				\Output{The array \textit{residual} containing all numbers $n\leq N$ satisfying the conditions  of $GGC_{m_1, m_2}$ for which there do not exist primes $p$ and $q$ such that $n=m_1p+m_2q$ and $m_1p\leq \alpha$.}
				
				\tcc{\textbf{Start Phase I: Unsegmented phase}}
				\tcc{Generating array 
					$primes$.} 
				\textbf{SmallPrimes($\max{(\lfloor \sqrt{\frac{N_{m_1, m_2}}{m_2}} \rfloor, \lfloor\frac{\alpha}{m_1}\rfloor )}$)}\;

				\tcc{Assigning values to arrays  $m_1p[r]$.}

				\textbf{Generatem1pr}($\alpha$)\;
				
				\tcc{Assigning values to array  $res$.}
				
				\textbf{GenerateResiduePattern}($m_1, m_2$)\;

				\tcc{\textbf{Start Phase II: Segmented phase}}
				
				\tcc{\textbf{Initialization}}

				Set arrays residual and $ism_2q$ empty\;
				
				$A\leftarrow 0$\;  

				\tcc{\textbf{Start segmented computation}}
				
				\While{$A<N$}{
					
					$B\leftarrow \min{\{A+\triangle, N\}}$\;

					\tcc{Saving the last $\alpha$ values of the boolean array $ism_2q$ generated in previous iteration in the boolean array $ism_2q\_old$.}	
					
					\If{$A>0$}{

						Save the last $\alpha$ elements of array ism2q in array ism$_2$q\_old\;  
						
					}
					
					\tcc{Assigning new values to array $ism_2q$.}
					\textbf{Generateism2q$(A,B)$}\;
					
					\tcc{Append array $ism_2q\_old$ to the beginning of array $ism\textsubscript{2}q$. In the updated array $ism_2q$, for every $0\leq i< \triangle+\alpha$, $ism_2q[i]=1$ if $A-\alpha +i$ is of the form $m_2q$ for some prime $q$ and $0$ otherwise.}	
					
					\If{$A>0$}{
						append\_before($ism_2q$; $ism_2q\_old$)\;
						
					}

					\tcc{Checking $GGC_{{}_{m_1,m_2}}$ in new interval.}
					
					\textbf{Check2$(A,B)$}\;
					
					$A\leftarrow A+\triangle$\;

				}		
			}
		\end{algorithm}
\end{small}}

\bibliographystyle{amsplain}


\end{document}



\section{}
\subsection{}

\begin{theorem}[Optional addition to theorem head]
\end{theorem}

\begin{proof}[Optional replacement proof heading]
\end{proof}

\begin{figure}
\includegraphics{filename}
\caption{text of caption}
\label{}
\end{figure}


\begin{equation}
\end{equation}

\begin{equation*}
\end{equation*}

\begin{align}
  &  \\
  &
\end{align}
